\documentclass[reqno]{amsart}
\usepackage{amsfonts}
\usepackage{mathrsfs}
\usepackage{amscd}
\usepackage{amsmath,amssymb,amsthm,latexsym,centernot}
\usepackage{enumerate}
\usepackage[english]{babel}
\usepackage{tikz-cd}
\usepackage{quiver}
\usepackage{enumitem}
\usepackage{stmaryrd} %For \owedge and \ovee symbols
\usepackage{hyperref}
\usepackage{color}

\newcommand{\cs}{\mathbf{RC}}
\newcommand{\csu}{\mathbf{RC}_{\mathbf{A}}}
\newcommand{\alg}{\mathbf{Alg}}
\newcommand{\calg}{\mathbf{cAlg}}
\newcommand{\rfd}{\mathbf{RFD}}
\newcommand{\frfd}{\mathbf{RFD}_{\mathbf{F}}}
\newcommand{\crfd}{\mathbf{cRFD}}

\newcommand{\cog}{\mathbf{Cog}}
\newcommand{\ccog}{\mathbf{cCog}}
\newcommand{\set}{\mathbf{Set}}
\newcommand{\vect}{\mathbf{Vect}}
\newcommand{\tops}{\mathbf{Top}}
\newcommand{\sob}{\mathbf{Sob}}

\newcommand{\cring}{\mathbf{cRing}}
\newcommand{\modu}{\mathbf{Mod}}
\newcommand{\modb}{\mathbf{Mod}_{\mathbf{B}}}
\newcommand{\cod}{\mathbf{Cod}}
\newcommand{\coh}{\mathbf{Coh}}
\newcommand{\tfc}{\mathbf{TFCoh}}

\newcommand{\sch}{\mathbf{Sch}}

\newcommand{\aff}{\mathbf{cAff}}

\newcommand{\mke}{\calg^\set}

\newcommand{\sq}{\square}
\newcommand{\n}{\mathfrak{n}}
\newcommand{\m}{\mathfrak{m}}
\newcommand{\p}{\mathfrak{p}}

\renewcommand{\O}{\mathcal{O}}
\newcommand{\M}{\mathfrak{M}}

\DeclareMathOperator{\spm}{Spm}

\DeclareMathOperator{\spc}{Spec}
\newcommand{\eps}{\varepsilon}
\newcommand{\dist}{\mathbf{T}}
\newcommand{\func}{\mathfrak{X}}
\newcommand{\emb}{\mathbf{T}}

\newcommand{\q}{\mathfrak{q}}
\newcommand{\lv}{\lVert}
\newcommand{\rv}{\rVert}
\newcommand{\dd}{\dag}
\newcommand{\Sec}{\mathcal{A}}
\newcommand{\F}{\mathfrak{F}}
\newcommand{\rcm}{\mathcal{F}}
\newcommand{\shm}{\mathscr{F}}

\newcommand{\til}{\widetilde}
\newcommand{\ha}{\widehat}

\DeclareMathOperator{\id}{id}

\DeclareMathOperator{\ann}{ann}

\newcommand{\invlim}{\varprojlim}
\newcommand{\dirlim}{\varinjlim}

\numberwithin{equation}{section}

\theoremstyle{plain}
\newtheorem{thm}[equation]{Theorem}
\newtheorem{cor}[equation]{Corollary}
\newtheorem{lem}[equation]{Lemma}
\newtheorem{prop}[equation]{Proposition}

\theoremstyle{definition}
\newtheorem{de}[equation]{Definition}

\newtheorem*{abst}{Abstract}

\newtheorem{exam}[equation]{Example}

\begin{document}

\title{A ringed-space-like structure on coalgebras for noncommutative algebraic geometry}

\author{So Nakamura}
%\author{Manuel L. Reyes}
%\address{Department of Mathematics\\ University of California, Irvine \\
%440M Rowland Hall\\ Irvine, CA 96267--3875 \\ USA}
%\email{mreyes57@uci.edu}
\date{\today}
%\thanks{This work was supported by NSF grant DMS-2201273}
\maketitle

\begin{abst}
Inspired by the perspective of Reyes' noncomutative spectral theory, we attempt to develop noncommutative algebraic geometry
by introducing \emph{ringed coalgebras}, which can be thought of as a noncommutative generalization of schemes over a field $k$. These objects arise from fully residually finite-dimensional(RFD) algebras introduced by Reyes and from schemes locally of finite type over $k$. The construction uses the Heyneman-Sweedler finite dual coalgebra and the Takeuchi underlying coalgebra. When $k$ is algebraically closed, the formation of ringed coalgebras gives a fully faithful functor out of the category of fully RFD algebras, as well as a fully faithful functor out of the category of schemes locally of finite type. The restrictions of these two functors to the category of (commutative) finitely generated algebras are isomorphic. Finally, we introduce \emph{modules over ringed coalgebras} and show that the category of finitely generated modules on a fully RFD algebra and that of coherent sheaves on a scheme locally of finite type can, if $k$ is algebraically closed, be fully faithfully embedded into the category of modules over the corresponding ringed coalgebras.

\end{abst} 
\setcounter{tocdepth}{2}
\makeatletter
\def\l@subsection{\@tocline{2}{0pt}{2.5pc}{5pc}{}}
\def\l@subsubsection{\@tocline{2}{0pt}{5pc}{7.5pc}{}} %Code to indent subsections for AMS packages
\makeatother
\tableofcontents

\setcounter{section}{-1}

\section{Introduction}
%The idea of viewing a ring as the set of functions over a space is originated in the study of $C^*$-algebras.

A duality between algebra and geometry appears in a various areas of mathematics. For example, there is a contravariant equivalence between the category of commutative (unital) $C^*$-algebras and the category of compact Hausdorff spaces, which is known as the Gelfand-Naimark duality. Similarly, the category of commutative reduced affine algebras over an algebraically closed field is contravariantly equivalent to the category of affine varieties over the field. Furthermore, the equivalence extends to the one betweeen commutative rings and affine schemes. For other examples, see Chapter 1 of \cite{K}.

The algebraic objects that appear in this context are all commutative rings. It would be natural to seek a noncommutative generalization of the correspondences; could we extend those fully faithful functors to some categories whose objects are noncommutative algebras? Unfortunately, the construction of the spectra of rings no longer work in noncommutative case, even for finite dimensional algebras. Indeed, the maximal spectra of the algebras $M_2(k)$ of $n\times n$ matrices over a field $k$ are empty. Moreover, \cite{R1} shows that any set-valued contravariant functor on the category of rings that extends the prime spectra functor on the category of commutative rings must send the square matrix algebras $M_n(\mathbb{C})$ to the empty set for all $n\geq3$. The maximum spectra functor on $C^*$-algebras has the same issue. The obstruction on the spectra functors would suggest that the underlying set of the``noncommutative spectrum" of commutative algebras should give us sets larger than the classical spectra of them. One could also observe that algebraic objects in the commutative cases play a role of rings of functions on the corresponding geometric objects; a commutative $C^*$-algebra is the algebra of continuous functions on its spectrum and a commutative affine $k$-algebra over an algebraically closed field $k$ is the algebra of regular functions on its maximal spectrum. For a noncommutative algebra, is it still possible to functorially construct an appropriate “spectrum" on which the elements of the algebra acts as "functions"? 

%Then the next question that should be asked would be what a noncommutative generalization of the correspondence should look like. The dualities mentioned at the beginning are given as categorical equivalences, and 

For a noncommutative $C^*$-algebra, one could consider the set of positive unital functionals on it, which is called the state space of the algebra. The state space is a convex subspace of the dual space of the algebra and contains the maximum spectrum of the algebra as its extremal points. As such, it gives a functor from the category of $C^*$-algebras to that of convex spaces and allows us to view a $C^*$-algebra as an algebra of functions on its state space. It is known that the functor can be turned into a fully faithful functor to an appropriate category using an additional structure on state spaces called rotations. For details, see \cite{AS}. 

%the state spaces with an additional structure called a rotation are constructed functorially and they give a fully faithful functor from the category of $C^*$-algebras to an appropriate category, 

%The construction of state spaces is functorial and 
%and allows us to view every The maximal spectrum of a commutative $C^*$-algebras can be recovered by taking the extremal points of the state space. 
%the geometrical interpretaion of algebras: viewing a noncommutative ring as a ring of functions on some kind of spectrum of it. 

On the other hand, an attempt to generalize algebraic geometry in a similar direction can be found in \cite{RV} and \cite{RV2}. The theory generalizes Nullstellensatz of commutative affine domains to that of affine PI algebras of fixed PI degree $n$. It leads a correspondence between affine PI algebras of PI degree $n$ and what are called $n$-varieties. However, the correspondence is not functorial in a usual sense: we need to restrict the morphisms of the category of PI algebras to the ones preserving surjective homomorphisms to the algebras of $n\times n$ matrices in order to establish a categorical equivalence between the algebras and the varieties.

%The algebraic objects that appear in the theory are finitely generated prime PI algebras of fixed PI degree $n$ over an algebraically closed field. The theory establishes a generalization of Nullstellensatz and as a result one obtain 
%an equivalence of categories between the category of such algebras with algebra homomorphism preserving serjectivity to the algebra of $n\times n$ matrices.

 %he category of such algebras and corresponding geometric objects, however, we need to restrict

%When the PI degree is 1, the theory is nothing but the usual affine algebraic geometry. 

%One assigns to such an algebra the topological space of algebra homomorphisms from $A$ to $M_n$, which makes the algebra an algebra of matrix-valued functions. 

%introduce $n$-varieties which generalize the algebraic varieties. 

%A topological space $U_{m, n}$ of $m$-tuple of $n\times n$ matrices that generate $M_n$ as an algebra plays a role of the affine $m$-space in usual algebraic geometry and $n$-varieties are defined to be the closed subsets of $U_{m, n}$ stable under PGL action.  However, the $n$-variety associated to a finitely generated prime PI algebra $A$ is . Because of this, 

In a recent paper \cite{R2} as well as in an article \cite[4.4]{R3}, where Reyes discusses noncommutative spectral theory, it is suggested that one should view coalgebras as ``noncommutative sets", in other words, a coalgebra is a candidate for the underlying structure of a ``noncommutative spectrum" of algebras. This perspective can be justified by the following reasons: first, the category $\set$ of sets can be fully faithfully embedded into $\cog$ by sending a set $X$ to the free vector space generated by $X$ with an appropriate coalgebra structure. Therefore coalgebras can be seen as a generalization of sets. Second, the underlying coalgebra construction introduced in \cite{T2} can be seen as a discretization of schemes in an appropriate sense(see Section 2.3 of \cite{R2}). Since the underlying coalgebras of affine schemes are the finite dual coalgebras of the corresponding commutative algebras, one should expect that the finite dual coalgebras of noncommutative algebras are discretizations of geometric objects corresponding to them. 

In this paper, we attempt to generalize algebraic geometry by using the finite dual coalgebra and the underlying coalgebra. Coalgebras are comonoid objects of the monoidal category of vector spaces. A finite dual coalgebra of an algebra as a vector space is a subspace of the dual space of the algebra. As such, every element of the algebra defines a linear map from the finite dual to the base field. The finite dual construction gives a contravariant functor from the category $\alg$ of algebras to the category $\cog$ of coalgebras. Following the philosophy of \cite{R2} and \cite[4.4]{R3}, we equip finite dual coalgebras with structures that play a role of closed subspaces and sheaves.

%The category $\set$ of sets can be fully faithfully embedded into the category $\cog$ so that coalgebras can be viewed as a generalization of sets. In this way, every algebra can be seen as an algebra of the functions on the finite dual coalgebra. 

%This paper is attempted to construct a theory that generalizes algebraic geometry using noncommutative RFD algebras over a fix field $k$. 

The outline of the paper is as follows: we first review basic notions such as the finite dual coalgebras of algebras and the underlying coalgebras of schemes in Section 1. We introduce a ringed-space-like structure on coalgebras and a category $\cs$ in Section 2 after discussing wedge products of subspaces of coalgebras. In section 3, we construct a fully faithful functor $(-)^\circ:\frfd^{op}\rightarrow\cs$ using the finite dual coalgebras where $\frfd$ is the full subcategory of $\alg$ whose objects are fully RFD algebras. In section 4 we construct a functor $\emb:\sch\rightarrow\cs$ using underlying coalgebras and show that it is compatible with $(-)^\circ$. Moreover, if the base field $k$ is algebraically closed, then the functor $\emb$ is fully faithful. The following is the main theorem(see \ref{com} and \ref{ffcom}):
\begin{thm}
The following diagram commutes up to natural isomorphism:
\[\begin{tikzcd}
	{\aff^{op}} && {\frfd^{op}} \\
	\sch^{lf} && \cs \\
	& \tops
	\arrow["incl.", hook, from=1-1, to=1-3]
	\arrow["\spc"', hook, from=1-1, to=2-1]
	\arrow["{(-)^\circ}", hook, from=1-3, to=2-3]
	\arrow["\emb", from=2-1, to=2-3]
	\arrow["{(-)(k)}"', from=2-1, to=3-2]
	\arrow["pts", from=2-3, to=3-2]
\end{tikzcd}\]
Here the round arrows stand for the fully faithfulness of functors. Furthremore, if $k$ is algebraically closed, then the functor $\emb$ is fully faithful.
\end{thm}
 In section 5, we introduce modules over ringed coalgebras and define a category ${}_C\modu$ of modules over $C$. We show that a category of finitely generated modules(resp. coherent sheaves) can be fully faithfully embedded to ${}_C\modu$ when $C=A^\circ$ for some RFD algebra $A$(resp.$C=\emb(X)$ for some scheme $X$ locally of finite type over an algebraically closed field $k$).

%These axioms are naturally induced when we generalize the construction of the Zariski topology on the maximum spectrum of an algebra to the finite dual coalgebra of it. One could consider the category $\cs$ whose objects are coalgebras equipped with a collection of subspaces and the morphisms are the coalgebra homomorphisms satisfying an axiom similar to that of continuous maps. The main features of this category are that $\cs$ contains both of the opposite category $\rfd^{op}$ of RFD algebras and the category $\tops$ of topological spaces as a reflective subcategory. The adjunction between $\cs$ and $\tops$ is compatible with that between $\cog$ and $\textbf{Set}$, which justifies to view $\cs$ the category of ``noncommutative topological spaces". In Section 3, we discuss analogous structure on comodules over some coalgebra and introduce the notion of ``continuous" Takeuchi equivalence. Finally, we show that Morita-equivalence is equivalent to ``continuous" Takeuchi equivalence.

\textbf{Acknowledgments.} The author is indebted to Manuel L. Reyes for his advice and unwavering guidance. He also thanks Thurman Ye for helpful comments on the draft.

%helpful advice and comments. He also thanks Thurman Ye for helpful comments on the draft.

\textbf{Notations.} Throughout this paper, we fix a field $k$. All vector spaces, algebras, coalgebras and tensor products $\otimes$ are the ones defined over $k$. All algebras have multiplicative identities and all algebra homomorphisms preserve them. Ideals always mean two-sided ideals. All modules(resp.comodules) over an algebra(resp.a coalgebra) are left modules(resp.right comodules) unless stated otherwise. For a vector space $V$, $V^*$ stands for the dual vector space of $V$, i.e., the vector space of all the linear maps from $V$ to $k$. We use the following notations for categories in this paper:
\begin{itemize}
\item $\set$: the category of sets and maps
\item $\tops$: the category of topological spaces and continuous maps
\item $\vect$: the cateory of vector spaces over $k$ and linear maps 
\item $\alg$: the category of algebras over $k$ and algebra homomorphisms
\item $\calg$: the full subcategory of $\alg$ whose objects are commutative algebras
%\item $\afd$: the full subcategory of $\calg$ whose objects are commutative finitely generated domains
\item ${}_A\modu$: the category of left $A$-modules and left $A$-module homomorphisms
\item $\modu_A$: the category of right $A$-modules and right $A$-module homomorphisms
\item $\cog$: the category of coalgebras and coalgebra homomorphisms
\item $\ccog$: the full subcategory of $\cog$ whose objects are cocommutative coalgebras
\item ${}_C\cod$: the category of left $C$-comodules and comodule homomorphisms
\item $\sch$: the category of schemes over $k$
\item $\sch^{lf}$: the full subcategory of $\sch$ whose objects are locally of finite type over $k$
\item $\coh_X$: the category of coherent sheaves over $X$
\end{itemize}

\section{Preliminaries}
In this section, we review some basic notions and constructions that we need for later sections. See \cite{DNR}, \cite{E}, \cite{S}, \cite{B} and \cite{T2} for details.

\subsection{The finite dual coalgebras}\label{cog0}
Recall that a triple $(C, \Delta, \varepsilon)$ is called a coalgebra if $C$ is a vector space, $\Delta:C\rightarrow C\otimes C$ and $\varepsilon: C\rightarrow k$ are linear maps that make the following diagrams commute:
\[\begin{tikzcd}
	C && {C\otimes C} && C && {C\otimes C} \\
	\\
	{C\otimes C} && {C\otimes C\otimes C} && {C\otimes C} && C
	\arrow["\Delta", from=1-1, to=1-3]
	\arrow["\Delta"', from=1-1, to=3-1]
	\arrow["{\Delta\otimes id_C}"', from=3-1, to=3-3]
	\arrow["{id_C\otimes \Delta}", from=1-3, to=3-3]
	\arrow["\Delta", from=1-5, to=1-7]
	\arrow["\Delta"', from=1-5, to=3-5]
	\arrow["{id_C\otimes\varepsilon}", from=1-7, to=3-7]
	\arrow["{\varepsilon\otimes id_C}"', from=3-5, to=3-7]
	\arrow["{id_C}"{description}, from=1-5, to=3-7]
\end{tikzcd}\]
A linear map $f:C\rightarrow C'$ between coalgebras $C$ and $C'$ is called a coalgebra homomorphism if it is compatible with the comultiplications and counits:
\[\begin{tikzcd}
	C && {C'} && C && {C'} \\
	\\
	{C\otimes C} && {C'\otimes C'} &&& k
	\arrow["f", from=1-1, to=1-3]
	\arrow["{\Delta_C}"', from=1-1, to=3-1]
	\arrow["{\Delta_{C'}}", from=1-3, to=3-3]
	\arrow["{f\otimes f}"', from=3-1, to=3-3]
	\arrow["f", from=1-5, to=1-7]
	\arrow["{\varepsilon_C}"', from=1-5, to=3-6]
	\arrow["{\varepsilon_{C'}}", from=1-7, to=3-6]
\end{tikzcd}\]
The category of coalgebras and coalgebra homomorphisms will be denoted by $\textbf{Cog}$.

If $C=(C, \Delta, \varepsilon)$ is a coalgebra, then the dual vector space $C^*$ is an algebra whose multiplication is defined by
$$(f\cdot g)(x):=(f\otimes g)(\Delta(x))$$
for all $x\in C$. The dual map of a coalgebra homomorphism is an algebra homomorphism and it gives a functor $(-)^*:\cog\rightarrow\alg^{op}$.

%A large class of examples of coalgebras
On the other hand, coalgebras can be produced from algebras. For an algebra $A$, one defines a subspace
$$A^\circ:=\{\phi\in A^*|\exists\ I\subset\ker\phi\ s.t.\ I\subset A\ \text{is a two-sided ideal with finite codimension}\}$$
of the dual vector space $A^*$. We can show that the dual map $m^*$ of the multiplication $m:A\otimes A\rightarrow A$ on the algebra $A$ induces a linear map $\Delta: A^\circ\rightarrow A^\circ\otimes A^\circ$ that makes $A^\circ$ a coalgebra together with $\varepsilon:A^\circ\rightarrow k$, the evaluation at $1\in A$. Note that if $\phi\in A^\circ$, then there exist collections of finitely many elements $\{\psi_{1, i}\}, \{\psi_{2, i}\}\subset A^\circ$ such that
$$\phi(ab)=\sum_i\psi_{1, i}(a)\psi_{2, i}(b)$$
for all $a, b\in A$. 

The coalgebra $A^\circ=(A^\circ, \Delta, \varepsilon)$ is called the finite dual coalgebra of $A$. It was introduced in \cite{HS}. The finite dual coalgebra has been considered as a candidate for a “noncommutative space" and studied in \cite{LB}, \cite{KS} and \cite{R2}.

The construction of the finite dual is functorial: if $\varphi:A\rightarrow B$ is an algebra homomorphism, then the dual linear map $\varphi^*:B^*\rightarrow A^*$ restricts to a coalgebra homomorphism $B^\circ\rightarrow A^\circ$. Therefore it induces a functor $(-)^\circ:\alg^{op}\rightarrow\cog$. This is a right adjoint of the functor $(-)^*:\cog\rightarrow\alg^{op}$, i.e., there is a natural bijection
$$\alg(A, C^*)\simeq\cog(C, A^\circ).$$

The unit $i_C: C\rightarrow C^{*\circ}$ of the adjunction is always injective whereas the counit $i_A:A\rightarrow A^{\circ*}$ is not necessarily. An algebra $A$ is called a \emph{residually finite-dimensional (RFD) algebra} if the cannonical linear map $i_A:A\rightarrow A^{\circ*}$ that sends every element in $A$ to its evaluation on $A^\circ$ is injective. This is equivalent to the condition that the intersection of all ideals of finite codimension in $A$ equals $0$, i.e., $A$ is a subdirect product of finite dimensional algebras. 

Any subalgebra of an RFD algebra is also RFD. Indeed, if $A$ is RFD algebra and $B\subset A$ is a subalgebra, then
$$\bigcap_{\dim(B/J)<\infty}J\subset\bigcap_{\dim(A/I)<\infty}I\cap B=0$$
where $I$ and $J$ are finite codimensional ideals of $A$ and $B$ respectively. The second condition implies that the product algebra $\prod A_i$ of any collection $\{A_i\}$ of RFD algebras is RFD. 

Denoting by $\rfd$ the full subcategory of $\alg$ whose objects are RFD algebras, we can restrict the adjunction explained above restricts to obtain the one between $\cog$ and $\textbf{RFD}^{op}$ since the dual algebra $C^*$ of a coalgebra $C$ is always RFD. We view the finite dual $A^\circ$ as a replacement of the maximal spectrum of $A$.

%Since the topology on $k$ is Hausdorff,  $V^*$ is always Hausdorff. In particular, $V^*$ is discrete if $\dim V<\infty$.The addition and scalar multiplication on $V^*$ are continuous in this topology and it becomes a topological vector space.

\subsection{Orthogonal spaces}\label{os}
In this subsection, we first review a correspondence of orthogonal spaces of $V$ and those of $V^*$ when the vector space $V$ is a coalgebra $C$. Then we consider the case $C=A^\circ$ for some RFD algebra $A$ and introduce some notations we will use in later sections. For details, see Chapter 1 of \cite{DNR}.

Let $C$ be a coalgebra. Recall that the finite topology on the dual vector space $C^*$ is the weakest topology that makes every evaluation $C^*\rightarrow k$ given by $\phi\mapsto\phi(x), x\in C$ continuous where $k$ has the discrete topology. The algebra structure on $C^*$ is indeed compatible with the finite topology on $C^*$ and hence $C^*$ is a topological algebra.

For every $T\subset C$ and $S\subset C^*$, we define the following subspaces of $C^*$ and $C$ respectively:
$$T^\perp:=\{\phi\in C^*|\phi(T)=0\}$$
and
\begin{eqnarray*}
S^\perp&:=&\{v\in C|\forall \phi\in S\ \phi(v)=0\}\\
&=&\bigcap_{\phi\in T}\ker\phi.
\end{eqnarray*}
These subspaces will be called the \emph{orthogonal spaces} of $T$ and $S$ respectively. The subspace $T^\perp\subset C^*$ is closed in the finite topology on $C$ and the assignments $T\mapsto T^\perp$ and $S\mapsto S^\perp$ give a bijective correspondence between subspaces of $C$ and closed subspaces of $C^*$. Finite dimensional subspaces of $C^*$ are closed.

A similar construction naturally appears in the fundamental of affine algebraic geometry. Suppose that $A$ is a commutative finitely generated algebra over an algebraically closed field $k$. We may identify the maximal spectrum $\spm(A)$, i.e., the set of maximal ideals of $A$, with the set of algebra homomorphisms from $A$ to $k$:
$$\spm(A)=\{\phi:A\rightarrow k|\phi\ \text{is an algebra homomorphism}\}\subset A^\circ.$$
A subset $T$ of the maximal spectrum $\spm(A)$ is Zariski closed if and only if $T$ is a vanishing set
$$Z(S):=\{\phi\in \spm(A)| \phi(S)=0\}=S^\perp\cap\spm(A)$$ 
of some subset $S\subset A$. Subsets of $\spm(A)$ of this form are called \emph{algebraic sets} and these subsets satisfy the axiom of closed subsets of topological spaces. Since $k$ is algebraically closed, the assignment $S\mapsto Z(S)$ gives a bijective correspondence between the radical ideals of $A$ and the closed subsets of $\spm(A)$, which is known as Hilbert's Nullstellensatz. Moreover, the set $R(A)$ of the radical ideals of $A$ and the set $C(\spm(A))$ of the closed subsets of $\spm(A)$ can be made into frames, i.e., lattices with finite meets and arbitrary joins. The order on $R(A)$ is given by the usual containment and the order on $C(\spm(A))$ is given by the opposite order of the usual containment. Under these lattice structures, the assignment $S\mapsto Z(S)$ does not only define a bijective map but also preserves finite meets and arbitrary joins, meaning that the two frames $R(A)$ and $C(\spm(A))$ are isomorphic.

Now we consider the case where $A$ is an RFD algebra and hence $A\subset A^{\circ*}$(here we do not assume that the base field $k$ is algebraically closed). We introduce some notations following the commutative affine case. A subspace
$$Z^\circ(S):=S^\perp=\{\phi\in A^\circ| \phi(S)=0\}$$ 
will be called a \emph{vanishing space} of a subset $S\subset A$ and a subspace in this form will be called an \emph{algebraic subspace}. We set 
$$\overline{S}:=S^{\perp\perp}\cap A=\bigcap_{\dim A/I<\infty} (S+I)$$
where $I$ runs over all the ideals of $A$ of finite codimension. The subspace $\overline{S}$ is the topological closure of $S$ in the relative topology on $A\subset A^{\circ*}$ of the finite topology on $A^{\circ*}$. It holds that $Z^\circ(\overline{S})=Z^\circ(S)$ for any subset $S\subset A$. A subspace $S\subset A$ will be called a closed subspace if $\overline{S}=S.$ If $\dim S<\infty$, then the subspace $S\subset A$ is closed(\cite[1.2.12]{DNR}). 

The assignment $S\mapsto Z^\circ(S)$ for closed subspaces $S\subset A$ defines a bijective correspondence between the closed subspaces of $A$ and the algebraic subspaces of $A^\circ$. Indeed, the inverse is given by $T\mapsto T^\perp=\bigcap_{\phi\in T}\ker\phi$ for algebraic subspaces $T\subset A^\circ$. If $S\subset A$ is a closed subspace, then we have
$$Z^\circ(S)^{\perp}=\bigcap_{\dim A/I<\infty} (S+I)=S.$$
On the other hand, if $T$ is an algebraic subspace of $A^\circ$, then $T=Z^\circ(S)$ for some subspace $S\subset A$. Thus we have
$$Z^\circ(T^{\perp})=Z^\circ\left(\bigcap_{\dim A/I<\infty} (S+I)\right)=Z^\circ(S)=T.$$
Furthermore, if $I\subset A$ is an ideal, then the algebraic subspace $Z^\circ(I)\subset A^\circ$ is a subcoalgebra(\cite[1.5.24]{DNR}). A subcoalgebra in this form will be called an \emph{algebraic subcoalgebra}.
The assignment $I\mapsto Z^\circ(I)$ gives a correspondence between the closed ideals of $A$ and algebraic subcoalgebras of $A^\circ$. This can be seen as a generalization of the Nullstellensatz.

%As a topological algebra, $C^*$ is a psudo-compact algebra under the finite topology and the ideals of finite codimension of $C^*$ form a system of neighborhood of $0$. 

%\subsection{Pointlike elements of coalgebras}

\subsection{Underlying coalgebras of schemes}\label{ucos}
We now briefly review the underlying coalgebras introduced in Chapter 2 of \cite{T2}. Recall that a coalgebra $C$ is said to be \emph{cocomutative} if the following diagram commutes:
\[\begin{tikzcd}
	& C \\
	{C\otimes C} && {C\otimes C}
	\arrow["\Delta"', from=1-2, to=2-1]
	\arrow["\Delta", from=1-2, to=2-3]
	\arrow["\tau"', from=2-1, to=2-3]
\end{tikzcd}\]
where $\tau$ is the twisting map defined by $\tau(x\otimes y)=y\otimes x$ for all $x, y\in C$. We denote by $\ccog$ the full subcategory of $\cog$ whose objects are cocomutative coalgebras. An algebra is commutative if and only if its finite dual $A^\circ$ is cocomutative. Likewise, $C$ is cocomutative if and only if its dual $C^*$ is commutative. In this way the adjunction between $\cog$ and $\alg$ naturally restricts to the one between $\ccog$ and $\calg$
where $\calg$ is a full subcategory of $\alg$ whose objects are commutative algerbas.

We say a coalgebra is \emph{pointed} if all simple subcoalgebras, i.e., the nonzero subcoalgebras that are minimal with respect to the containment, are of 1-dimensional. It is \emph{irreducible} if any two nonzero subcoalgebras intersect nontrivially. A cocomutative coalgebras can always be written as a direct sum of its pointed irreducible subcoalgebras(\cite[8.0.7]{T2}).

Every scheme $X$ over $k$ defines a functor $X:\calg\rightarrow\set$ by setting
$$X(A):=\sch(\spc(A), X).$$
This assignment naturally defines a fully-faithful functor from the category $\sch$ of ($k$-)schemes to the functor category $\mke:=\textbf{Func}(\calg, \set)$, which allows us to view $\sch$ as a full subcategory of $\mke$. The underlying coalgebra $\dist(\func)$ of a functor $\func:\calg\rightarrow\set$ is defined to be a cocomutative coalgebra satisfying
$$\ccog(C, \dist(\func))\simeq \func(C^*)$$
for all cocomutative coalgebras $C$ of finite dimension. Note that such $\dist(\func)$ is unique up to isomorphism if it exists. If $\func=X$ for some scheme $X$ over $k$, the coalgebra $\dist(X)$ exists and is given by
$$\dist(X)\simeq\bigoplus_{x\in\lv X\rv} \O_{X, x}^\circ$$
where 
$$\lv X\rv=\{x\in X|[\kappa(x):k]<\infty\}.$$
Here $\kappa(x)$ stands for the function field at $x$. In the above isomorphism, each $\O_{X, x}^\circ$ is a pointed irreducible coalgebra and any pointed irreducible subcoalgebra of $\dist(X)$ is in this form. The coalgebra $\O_{X, x}^\circ$ is called the \emph{tangent coalgebra} of $X$ at $x$.

If $X$ is affine and $X\simeq\spc(A)$ for some commutative RFD algebra $A$, then 
$$\dist(\spc(A))\simeq\bigoplus_{\m\in \lv\spc(A)\rv} A^\circ_\m\simeq A^\circ.$$

The construction of $\dist$ defines a functor $\dist:\sch\rightarrow\cog$. The image of the objects of $\sch$ under the functor $\dist$ lies in $\ccog$. If $f:X\rightarrow Y$ is a morphism of schemes, then $\dist(f):\dist(X)\rightarrow \dist(Y)$ is the direct sum of the finite duals $f_x^\circ:\O_{X, x}^\circ\rightarrow\O_{Y, f(x)}^\circ$ followed by the inclusion to $\dist(Y)$. This functor commutes with finite limits by the definition of $\dist$ since the subcategory $\sch$ of $\mke$ is closed under finite limits. We point out that the underlying coalgebras of schemes are viewed in \cite{R2} as  “underlying discrete objects" of schemes. 

%It is shown in a part of 2.2.5 of \cite{T2} that if $X$ and $Y$ are scheme and $f:X\rightarrow Y$ is a monomorphism of scheme, then $\dist(f):\dist(X)\rightarrow\dist(Y)$ is injective. 2.2.7 shows that the converse is also true if $X$ and $Y$ are locally of finite type. Since the full subcategory $\sch$ of $\mke$ is closed under limits, so it is equivalent to $f$ being monic in $\sch$. 

Finally, we close this subsection with a discussion of some properties of the functor $\dist:\sch\rightarrow\cog$. The following diagram commutes up to natural isomorphism(\cite{T2}, 2.1.6, 2.1.7):
\[\begin{tikzcd}
	{\crfd^{op}} & {\rfd^{op}} \\
	\sch & \cog
	\arrow["incl.", hook, from=1-1, to=1-2]
	\arrow["\spc", hook, from=1-1, to=2-1]
	\arrow["(-)^\circ", from=1-2, to=2-2]
	\arrow["\dist"', from=2-1, to=2-2]
\end{tikzcd}\]
where $\crfd$ is the full subcategory of $\rfd$ whose objects are commutative RFD algebras.

The functor $\rfd^{op}\rightarrow\cog$ is faithful but is not full. The faithfulness follows from the commutativity of the following diagram where $A$ and $B$ are RFD algebras and $\varphi:A\rightarrow B$ is an algebra homomorphism:
\[\begin{tikzcd}
	{A^{\circ*}} & {B^{\circ*}} \\
	A & B
	\arrow["{\varphi^{\circ*}}", from=1-1, to=1-2]
	\arrow[hook, from=2-1, to=1-1]
	\arrow["\varphi"', from=2-1, to=2-2]
	\arrow[hook, from=2-2, to=1-2]
\end{tikzcd}\]
However, we have
$$\rfd(k[x], A^{\circ*})\simeq\cog(A^\circ, k[x]^\circ)$$
by adjunction and the left hand side(LHS) is strictly larger than $\rfd(k[x], A)$ as a set if $A$ is of infinite dimension. 

Proposition 2.2.7 of \cite{T2} shows that the restriction of the functor $\dist$ to $\sch^{lf}$, the full subcategory of $\sch$ whose objects are locally of finite type, is faithful, but it is still not full. Again, $\sch(\spc(A), \spc(k[x]))=\rfd(k[x], A)$ and this is in general a proper subset of the hom-sets
$$\cog(\dist(\spc(A)), \dist(\spc(k[x])))=\cog(A^\circ, k[x]^\circ).$$

\section{Category $\cs$ of ringed coalgebras}

We start this section with a discussion of some properties of wedge products of subspaces of coalgebras. The wedge product is an operation on the sets of the subspaces of a coalgebra and can be seen as an analogue of the union operation on the powerset of a set. We discuss relationships with wedge products and some set-theoretic operations such as intersections and preimages of subspaces. As usual, $\cog$ denotes the category of coalgebras and coalgebra homomorphisms and $\sch$ denotes the category of $k$-schemes and scheme morphisms.

%Viewing coalgebras as a noncommutative generalization of sets, we define a topology-like structure on coalgebras using the wedge products.

%We especially look into set-theoretic properties of the wedge products of subspaces of coalgberas since they replace the unions of subsets when we view coalgebras as a noncommutative generalization of sets. We also discuss the construction of finite dual comodule in a full detail. It is briefly mentioned in a full detail so it is discussed in this section. %and so on and show some properties of them that we use later on. 

\subsection{Wedge products as noncommutative unions}
For two subspaces $S_1, S_2$ of an algebra $A$, one can consider the product space $S_1S_2$ of them:
$$S_1S_2:=\left\{\sum_{1\leq i\leq n}s_is'_i\in A\mid\exists n\geq 1, s_i\in S_1, s'_i\in S_2\right\}.$$
 The dualized operation on the subspaces of a coalgebra is the wedge product. Let $C=(C, \Delta, \eps)$ be a coalgebra and $V, W\subset C$ be a subspace. Then the wedge product is defined by
$$V\vee_C W:=\Delta^{-1}(V\otimes C+C\otimes W).$$
%The usual notation of this notion in the literature is $\wedge$. We use the notation $\vee$ in order to avoid any confusion due to using $\wedge$ and $\cap$ because we also consider the intersection of subspaces of a coalgebra at the same time in this paper. 

This subspace $V\vee W$ can be verified to be equal to $(V^\perp W^\perp)^\perp$ where the product $V^\perp W^\perp$ is taken as subspaces of $C^*$. In particular, it is associative, i.e., 
$$(V_1\vee V_2)\vee V_3=V_1\vee(V_2\vee V_3).$$
If $V$ and $W$ are subcoalgebras, then so is $V\vee W$.

%Unlike union of subsets of a set, the wedge product is neither commutative nor idempotent in general.

%Recall that the free vector space $kX$ has a coalgebra structure given by the comultiplication and the counit. For a set $X$, the vector space $kX$ formally generated by $X$ is a coalgebra with the comultiplication defined by $x\mapsto x\otimes x$ and the counit defined by $x\mapsto 1$. This yields a fully faithful functor $\set\hookrightarrow\cog$. 

The wedge product of subspaces can be seen as noncommutative generalization of union of the subsets of a set. We recall a coalgebra structre on a free vector space generated by a set. Let $X$ be a set and $kX$ be a vector space generated by $X$. Then the map $kX\rightarrow kX\otimes kX$ induced by $x\mapsto x\otimes x$ and the map $kX\rightarrow k$ induced by $x\mapsto 1$ make $kX$ a coalgebra. It can be easily seen that this defines a functor $k[-]:\set\rightarrow\cog$. It has a right adjoint given by taking the group-like elements of coalgebras. Recall that a nonzero element $x\in C\backslash\{0\}$ is called a group like element if $\Delta(x)=x\otimes x$. The subset of group like elements are denoted by $pts(C)$. Since coalgebra homomorphisms take group-like elements of domains to those of codomains, it defines a functor $pts:\cog\rightarrow\set$.
$$\cog(kX, C)\simeq\set(X, pts(C)).$$
Moreover, $pts(kX)=X$ for any set $X$ so that 
$$\cog(kX, kY)\simeq\set(X, pts(kY))=\set(X, Y).$$
Thus $k[-]$ is fully-faithful and $\set$ is a reflective subcategory of $\cog$. In this way, coalgebras can be viewed as noncommutative generalization of sets and its duals can be viewed as the algebras of $k$-valued functions on them. The functor $pts$ is compatible with the functor $(-)(k):=\sch(\spc(k), -):\sch\rightarrow\set$ that gives the $k$-rational points of $X$, i.e., the following commutes up to natural isomorphism:
\[\begin{tikzcd}
	\sch && \cog \\
	\\
	& \set
	\arrow["\dist", from=1-1, to=1-3]
	\arrow["(-)(k)"', from=1-1, to=3-2]
	\arrow["pts", from=1-3, to=3-2]
\end{tikzcd}\]

The subsets of $X$ bijectively correspond to the subcoalgebras of $kX$ under the assignment $S\mapsto kS$. Moreover, it holds that 
$$kS\vee kT=k[S\cup T]$$
for any subsets $S, T\subset X$. This gives an isomorphism of semi-lattices between the powerset of $X$ with union and the set of subcoalgebras of $kX$ with wedge product. 

The above equation implies that
$$(kS\vee kT)\cap X=S\cup T$$
and this can be proved in a general case as follows: 
\begin{lem}
Let $C$ be a coalgebra and $V_1, V_2\subset C$ be subspaces. Then it holds that $$(V_1\vee V_2)\cap pts(C)=(V_1\cap pts(C))\cup (V_2\cap pts(C)).$$
\end{lem}
\begin{proof}
It is immediate to see that the right hand side (RHS) contains the left hand side (LHS). Let $x$ be an element of the LHS and assume that $x\notin V_1$. Then there exists $\phi\in C^*$ such that $\phi(V_1)=0$ and $\phi(x)\neq0$. Let $\psi\in C^*$ be an element such that $\psi(V_2)=0$. Then we have $\phi(x)\psi(x)=0$ since $x\in (V_1^\perp V_2^\perp)^\perp$. Therefore $\psi(x)=0$. This shows that $x\in V_2^{\perp\perp}=V_2$. The symmetric argument shows that $x\in V_1$ if $x\notin V_2$. Thus $x$ is an element of the RHS.
\end{proof}
%\section{Coalgebras as noncommutative sets}

The rest of this section will be spent for studying properties of wedge products. We first show that the wedge products distribute over intersection(See\cite[2.4]{E}).

\begin{lem}\label{lin}
Let $V$ and $W$ be vector spaces.
\begin{enumerate}
\item  Let $V'\subset V$ be a subspace and $x=\sum_s v_s\otimes w_s$ be an element of $V\otimes W$ where $v_s\in V$ and $w_s\in W\backslash\{0\}$. If $x\in V'\otimes W$, then $v_s\in V'$ for all $s$.
\item Let $\{V_i\}$ be a collection of subspaces of $V$. Then we have
$$\bigcap_i (V_i\otimes W)=(\bigcap_i V_i)\otimes W$$
as subspaces of $V\otimes W$.
\end{enumerate}
\end{lem}
\begin{proof}
(1) Let $B_1=\{u_j\}$ be a basis of $V'$ and $B_2=\{v_k\}\subset V$ be a subset disjoint from $B_1$ such that $B_1\cup B_2$ is a basis of $V$. Let $B_3=\{w_p\}$ be a basis of $W$. Then 
$$v_s=\sum_j c^s_j u_j+\sum_k c^s_{k} v_k,\quad w_s=\sum_p d^s_p w_p$$
where $c^s_j, c^s_k, d^s_p\in k$. %and $u^s_j\in B_1, v^s_k\in B_2$ and $w^s_p\in B_3$. 
Here we may assume that $d^s_p\neq0$ for all $s$ and $p$. Therefore
\begin{eqnarray*}
x&=&\sum_s v_s\otimes w_s\\
&=&\sum_{s, j, p}c^s_jd^s_pu_j\otimes w_p+\sum_{s, k, p}c^s_kd^s_pv_k\otimes w_p
\end{eqnarray*}
Since $x$ and the first term in the above expression are elements of $V'\otimes W$, $c^s_kd^s_p=0$ for all $s, k$ and $p$. Since $d^s_p\neq0$, we have that $c^s_k=0$ for all $s$ and $k$. Thus
$$v_s=\sum_j c^s_j v_j+\sum_k c^s_{k} v_k=\sum_j c^s_j v_j\in V'$$
for all $s$.\\
(2) It is enough to show that the LHS is contained in the RHS. Let $x=\sum_sv_s\otimes w_s$ be an element of the LHS where $v_s\in V$ and $w_s\in W$. We may assume that $w_s\neq0$ for all $s$. For each $i$, we know that $x=\sum_sv_s\otimes w_s\in V_i\otimes W$. By (1), we have $v_s\in V_i$ for all $s$. Therefore $v_s\in\bigcap_i V_i$ for all $s$ and the proof completes.
\end{proof}

\begin{lem}
Let $C$ be a coalgebra, $\{V_i\}$ be a collection of subspaces of $C$ and $W\subset C$ be a subspace. Then it holds that
$$(\bigcap_i V_i)\vee_C W=\bigcap_i (V_i\vee_C W)\ \text{and}\ W\vee_C(\bigcap_i V_i)=\bigcap_i (W\vee_C V_i).$$
\end{lem}
\begin{proof}
We only show the first equality 
\begin{eqnarray*}
(\bigcap_i V_i)\vee_C W=\bigcap_i (V_i\vee_C W)
\end{eqnarray*}
since the argument for the second equality is similar. It is easy to see that the right hand side (RHS) contains the left hand side (LHS). To see the containment of the other direction, we consider the following commutative diagram:
\[\begin{tikzcd}
	C && {(C/\bigcap_i V_i)\otimes (C/W)} \\
	\\
	&& {\prod_i(C/V_i\otimes C/W)}
	\arrow["\Delta", from=1-1, to=1-3]
	\arrow["{\prod\circ\Delta}"', from=1-1, to=3-3]
	\arrow["\prod", from=1-3, to=3-3]
\end{tikzcd}\]
Here the map $\Delta$ is the composite of the comultiplication of $C$ and the natural quotient map $C\otimes C\twoheadrightarrow (C/\bigcap_i V_i)\otimes C/W$.  The map $\prod$ is the product of the natural quotient maps $(C/\bigcap_i V_i)\otimes C/W\twoheadrightarrow (C/V_i)\otimes (C/W)$. By definition, the wedge product
$$(\bigcap_i V_i)\vee_C W$$
is the kernel of $\Delta$ and the intersection
$$\bigcap_i (V_i\vee_C W)$$
is the kernel of $\prod\circ\Delta$. Thus it suffices to show that the canonical map 
$$\Pi:(C/\bigcap_i V_i)\otimes C/W\rightarrow \prod_i(C/V_i\otimes C/W)$$
is injective. The kernel of this map is 
$$\bigcap_j ((V_j/\bigcap_i V_i)\otimes C/W)$$
where the intersection is taken in $(C/\bigcap_i V_i)\otimes C/W$.
By (2) of lemma \ref{lin} we have 
$$\bigcap_j ((V_j/\bigcap_i V_i)\otimes C/W)=(\bigcap_j (V_j/\bigcap_i V_i))\otimes (C/W)=0\otimes C/W=0.$$
Thus, the map $\Pi$ is injective. This implies that the $\bigcap_i (V_i\vee_C W)$ contains $(\bigcap_i V_i)\vee_C W$. The proof is complete.
\end{proof}

\begin{cor}\label{AQ}(cf. \cite{E} Exercise 2.5.21)
Let $A$ be an RFD algebra.
\begin{enumerate}
\item If $\{S_i\}$ is a collection of linear subspaces of $A$, then $\bigcap _iZ^\circ(S_i)=Z^\circ\left(\sum_i S_i\right)$.
\item If $S_1, S_2\subset A$ are linear subspaces, then $Z^\circ(S_1)\vee Z^\circ(S_2)=Z^\circ(S_1S_2)$.
\end{enumerate}
\end{cor}
\begin{proof}
(1) Immediate.\\
(2) By the previous lemma, we have
\begin{eqnarray*}
Z^\circ(S_1)\vee Z^\circ(S_2)&=&\left(\bigcap_{s\in S_1} Z^\circ(s)\right)\vee\left(\bigcap_{s'\in S_2} Z^\circ(s')\right)=\bigcap _{s\in S_1, s'\in S_2}(Z^\circ(s)\vee Z^\circ(s')).
\end{eqnarray*}
Here 
$$Z^\circ(s)\vee Z^\circ(s')=(Z^\circ(s)^\perp Z^\circ(s')^\perp)^\perp=(kss')^\perp=Z^\circ(ss').$$
The first equality follows because $ks, ks'\subset A$ are of finite dimension and hence closed. Therefore
\begin{eqnarray*}
Z^\circ(S_1)\vee Z^\circ(S_2)=\bigcap _{s\in S_1, s'\in S_2}Z^\circ(ss')=Z^\circ(\{\sum_i s_is'_i\mid s_i\in S_1, s'_i\in S_2\})=Z^\circ(S_1S_2)
\end{eqnarray*}
\end{proof}

Therefore algebraic subspaces of $A^\circ$ is closed under arbitrary intersection and finite wedge products. In particular they form a quantale, i.e., a monoid object in the category of complete upper semi-lattices. The assignment $S\mapsto Z^\circ(S)$ is not only a bijection but indeed an isomorphism of quantales. Here the quantale structure on the set of the closed subspaces of $A$ is given by
$$\bigvee S_i=\overline{\sum_i S_i},\ S_1*S_2:=\overline{S_1S_2}$$
where $\overline{(-)}$ is the closure in the finite topology on $A\subset A^{\circ*}$. This correspondence can be seen as a noncommutative version of the Nullstellensatz.

Next, we consider the relationship between wedge products and  coalgebra homomorphisms. Note that taking image is not compatible in general: it holds that
$$f(V\vee W)\subset f(V)\vee f(W)$$
but the containment is not necessarily an equality. On the other hand, preimage under a coalgebra homomorphism is compatible with wedge product.

\begin{lem}
Let $C, D$ be coalgebras and $f:C\rightarrow D$ be a coalgebra homomorphism. Then for all subspaces $V\subset D$ and $W\subset D$, it holds that
$$f^{-1}(V\vee_D W)=f^{-1}(V)\vee_C f^{-1}(W).$$
\end{lem}
\begin{proof}
Note that $V\vee_D W$ is the kernel of the composition $D\rightarrow D\otimes D\rightarrow D/V\otimes D/W$ of the comultiplication and the natural quotient linear map. Consider the following commutative diagram whose rows are both exact sequences:
\[\begin{tikzcd}
	0 & {V\vee_D W} & D & {D/V\otimes D/W} & 0 \\
	0 & {f^{-1}(V)\vee_C f^{-1}(W)} & C & {C/f^{-1}(V)\otimes C/f^{-1}(W)} & 0 \\
	{f^{-1}(V\vee_D W)}
	\arrow[from=1-1, to=1-2]
	\arrow[hook, from=1-2, to=1-3]
	\arrow["\Delta_D", two heads, from=1-3, to=1-4]
	\arrow[from=1-4, to=1-5]
	\arrow[from=2-1, to=2-2]
	\arrow["f", from=2-2, to=1-2]
	\arrow[hook, from=2-2, to=2-3]
	\arrow["f"', from=2-3, to=1-3]
	\arrow["\Delta_C", two heads, from=2-3, to=2-4]
	\arrow["{f\otimes f}"', hook, from=2-4, to=1-4]
	\arrow[from=2-4, to=2-5]
	\arrow["f", shift left, from=3-1, to=1-2]
	\arrow[shift right=2, hook, from=3-1, to=2-3]
\end{tikzcd}\]
Let $x\in f^{-1}(V\vee_D W)$. Then $\Delta_D(f(x))=0$. By the commutativity, we have $f\otimes f\circ \Delta_C(x)=0$. Since $f\otimes f$ is injective, $\Delta_C(x)=0$ so that $x\in \ker\Delta_C=f^{-1}(V)\vee_C f^{-1}(W)$. Thus $f^{-1}(V\vee_D W)\subset f^{-1}(V)\vee_C f^{-1}(W).$ If $x\in f^{-1}(V)\vee_C f^{-1}(W)$, then $x$ is mapped to $V\vee_D W$ under $f$ so that it must be contained in $f^{-1}(V\vee_D W)$. Therefore $f^{-1}(V\vee_D W)=f^{-1}(V)\vee_C f^{-1}(W).$
\end{proof}

\begin{cor}\label{int/wed}
Let $C$ be a coalgebra, $V_1, V_2$ be subspaces, and $D$ be subcoalgebras of $C$. Then it holds that
$$(V_1\vee_C V_2)\cap D=(V_1\cap D)\vee_D(V_2\cap D)$$
\end{cor}
\begin{proof}
Apply the preceding lemma to the inclusion $D\hookrightarrow C$ and subspaces $V_1, V_2\subset C$.
\end{proof}

%In this section we define a topology-like structure on a coalgebra and consider category theoretic relationships with topological spaces and RFD algebras. 

\begin{lem}\label{sumprod}
Let $\{C_i\}$ be a collection of coalgebras, $C=\bigoplus_i C_i$ be the direct sum and $V_i, V'_i\subset C_i$ be subspaces for all $i$.
Then 
$$\left(\bigoplus_i V_i\right)\vee_C\left(\bigoplus_i V'_i\right)=\bigoplus_i (V_i\vee_{C_i} V_i').$$
\end{lem}
\begin{proof}
The statement follows from the definition of wedge products and the equality
$$\left(\bigoplus V_i\right)\otimes C+C\otimes \left(\bigoplus V'_i\right)
=\bigoplus_i(V_i\otimes C+C\otimes V'_i).$$
\end{proof}

We close this section by discussing category-theoretic pullbacks in $\cog$. If $f:C_1\rightarrow C_2$ is a homomorphism of coalgebras and $D\subset C_2$ is a subcoalgebra, it is not true in general that the preimage $f^{-1}(D)\subset C_1$ is a subcoalgebra. For instance, $\ker f=f^{-1}(0)\subset C_1$ is not a subcoalgebra even though $0\subset C_2$ is unless $f$ is injective. This is because $\ker f\subset C_1$ is a coideal and any coideal that is a subcoalgebra is always $0$.

The category-theoretic pullback of a subcoalgebra under a coalgebra homomorphism has the following description:

\begin{lem}\label{pb}
Let $f:C_1\rightarrow C_2$ be a coalgebra homomorphism and $D\subset C_2$ be a subcoalgebra. Then the sum
$$f^\dd(D):=\sum_{D'\subset f^{-1}(D)}D'\subset C_1$$
where $D'$ runs all the subcoalgebras of $C_1$ contained in $f^{-1}(D)$ is a pullback of $D$ under $f$ in the category $\cog$. Furthermore, $f^\dd(D\cap D')=f^\dd(D)\cap f^\dd(D')$ for any subcoalgebras $D, D'\subset C_2$.
\end{lem}
\begin{proof}
Let $f':f^\dd(D)\rightarrow D$ be the restriction of $f$ to $f^\dd(D)$, $i:D\hookrightarrow C_2$ and $i':f^\dd(D)\hookrightarrow C_1$ the inclusions. Suppose that $D'$ is a coalgebra, $g_1:D'\rightarrow C_1$ and $g_2:D'\rightarrow D$ are coalgebra homomorphisms such that $f\circ g_1=i\circ g_2$.
 \[\begin{tikzcd}
	& {C_1} & {C_2} \\
	& {f^\dd(D)} & D \\
	{D'}
	\arrow["f", from=1-2, to=1-3]
	\arrow["i'", hook, from=2-2, to=1-2]
	\arrow["f'", from=2-2, to=2-3]
	\arrow["i"', hook, from=2-3, to=1-3]
	\arrow["{g_1}", from=3-1, to=1-2]
	\arrow["j", dotted, from=3-1, to=2-2]
	\arrow["{g_2}"', from=3-1, to=2-3]
\end{tikzcd}\]
Then $g_1(D')\subset f^{-1}(D)$. Since $g_1(D')$ is a subcoalgebra of $C_2$, it is contained in the subcoalgebra $f^\dd(D)\subset C_1$. Thus $g_1$ is corestricted to a coalgebra homomorphism $j:D'\rightarrow f^\dd(D)$ and it makes the diagram commute.

The second assertion holds because an intersection $D\cap D'$ is a pullback of the following commuting diagram of inclusions:
\[\begin{tikzcd}
	D & {C_2} \\
	{D\cap D'} & {D'}
	\arrow[hook, from=1-1, to=1-2]
	\arrow[hook, from=2-1, to=1-1]
	\arrow[hook, from=2-1, to=2-2]
	\arrow[hook, from=2-2, to=1-2]
\end{tikzcd}\]
\end{proof}

Recall that a subcoalgebra $C$ of a finite dual coalgebra $A^\circ$ is called algebraic if $C=Z^\circ(I)$ for some ideal $I\subset A$. If $C_1$ and $C_2$ are finite dual coalgebras of some algebras and $f$ is a finite dual of an algebra homomorphism, then the category-theoretic pullback of any algebraic subcoalgebra under $f$ is algebraic.  

\begin{lem}\label{algpb}
Let $\varphi:A\rightarrow B$ be an algebra homomorphism and $f:B^\circ\rightarrow A^\circ$ be its finite dual. Let $D=Z^\circ(I)\subset A^\circ$ be a subcoalgebra defined by an ideal $I\subset A$. Then $f^\dd(D)=Z^\circ(B\varphi(I)B)$. 
\end{lem}
\begin{proof}
Let $C\subset B^\circ$ be a subcoalgebra contained in $Z^\circ(\varphi(I))$ and $\phi\in C$. Using the Sweedler notation, write
$$(\Delta\otimes id)\circ\Delta(\phi)=\sum \psi_{(1)}\otimes\psi_{(2)}\otimes \psi_{(3)}$$ 
where $\psi_{(1)}, \psi_{(2)}$ and $\psi_{(3)}$ are elements in $C$. Note that 
$$\phi(xyz)=\sum \psi_{(1)}(x)\psi_{(2)}(y)\psi_{(3)}(z)$$
for all $x, y, z\in B$ since $\phi\in B^\circ$. In particular, we have
$$\phi(b_1\varphi(a)b_2)=\sum\psi_{(1)}(b_1)\psi_{(2)}(\varphi(a))\psi_{(3)}(b_2)=0$$
for any $a\in I, b_1, b_2\in B$ since $\psi_{i, j}\in C\subset Z^\circ(\varphi(I))$. Thus $\phi\in Z^\circ(B\varphi(I)B)$ and hence $C\subset Z^\circ(B\varphi(I)B)$.
\end{proof}
For two subcoalgebras $D, D'\subset C_2$, it is not true in general that $f^\dd(D\vee D')=f^\dd(D)\vee f^\dd(D')$. Indeed, let $A:=k[x]$ be the algebra of polynomials with a single variable, $B=M_2(k)$ be the matrix algebra and $\varphi:k[x]\rightarrow M_2(k)$ be an algebra homomorphism such that $\varphi(x)$ is a nonzero matrix whose square is the zero matrix. Take the finite dual $f:=\varphi^\circ:M^2(k)\rightarrow k[x]^\circ$ and the ideal $I:=(x)\subset k[x]$ generated by $x$. Note that $\varphi(I)\neq0$ but $\varphi(I^2)=0$. Then it follows by Lemma \ref{algpb} that
$$f^\dd(Z^\circ(I))=Z^\circ(B\varphi(I)B)=Z^\circ(M_2(k))=0$$
so that $f^\dd(Z^\circ(I))\vee f^\dd(Z^\circ(I))=0\vee0=0$.
However 
$$f^\dd(Z^\circ(I)\vee Z^\circ(I))=f^\dd(Z^\circ(I^2))=Z^\circ(B\varphi(I^2)B)=Z^\circ(0)=M^2(k).$$
This example also shows that a preimage of a subcoalgebra under a coalgebra homomorphism is not necessarily a subcaogebra; the image $Im(\varphi)\subset M_2(k)$ is of dimension one and hence $\ker f=Im(\varphi)^\perp$ is of dimension three. However, the subcoalgebras of $M^2(k)$ are $0$ or $M^2(k)$ so that $\ker f\subset M^2(k)$ cannot be a subcoalgebra. 

%write $M^2(k)=\ker f\oplus N$ for some subspace $N\subset M^2(k)$. There exists an element $A, B\in M_2(k)$ such that $A$ is $0$ on $\ker f$ and 

%the kernel $\ker f$ is nonzero since the image $Im(\varphi)\subset M_2(k)$ has a dimension one.

\subsection{Definition of $\cs$}

%The quantale $P_A$ of an RFD algebra $A$ is analogous to Mulvay's quantale Max($\mathcal{A}$), the set of all closed subspaces of $\mathcal{A}$, of a $C^*$-algebra $\mathcal{A}$. The quantale Max($\mathcal{A}$) can be viewed as a noncommutative generalization of compact Hausdorff spaces as 

% We also have $\ker\varepsilon=Z^\circ(k)$ and $A^\circ=Z^\circ(0)$ by definition. These properties can be seen analogously to the axioms of closed subsets of topological spaces. Moreover, the assignment $S\mapsto Z^\circ(S)$ gives a bijective map $P_A$ to $Q_A$ and the properties mentioned above means that they are isomorphic as monoid objects in the category of semi-lattices. 

%The difference is that Max($A$) is complete as a lattice whereas $P_A$ is not. We will not consider $C^*$-algebras in this paper, but it would be interesting problem to thinks about, for example, whether or not the assignment $A\mapsto P_A$ classifies the isomorphism classes of $C^*$-algebras.

In this section, we define ringed coalgebras, which can be thought of as noncommutative version of ringed spaces.

Recall the construction of affine varieties over an algebraically closed field $k$. If $A$ is a commutative reduced affine algebra, then $\spm(A)$ is endowed with a Zariski topology where closed subsets are of the form $Z(I)$ for some ideal $I\subset A$. Every element of $A$ defines a continuous map from $\spm(A)$ to $k$ and $A$ can be viewed as an algebra of $k$-valued functions on $\spm(A)$. Furthermore, the structure sheaf on $\spm(A)$ is defined using the localizations of $A$. 

We attempt to follow this process when $A$ is an RFD algebra. The philosophy of this paper is to view $A$ as an algebra of functions on the finite dual $A^\circ$. The closed subsets $Z(I)\subset\spm(A)$ for some ideal $I\subset A$ are replaced by subcoalgebras of the form $Z^\circ(I)$. However, there does not seem to be an obvious counterpart of open subsets of $\spm(A)$. Therefore we assign to every subcoalgebra $C\subset A^\circ$ a subalgebra $A_C\subset C^*$ in such a way that the assignment gives a functor $P_C^{op}\rightarrow\rfd$ where $P_C$ is the category of subcoalgebras of $C$ whose morphisms are the inclusions. 

%the pair $(A^\circ, Q_A)$ as an object analogous to a topological space. Roughly speaking, the finite dual coalgebra $A^\circ$ is the underlying object of this pair and $Q_A$ plays a role of “topological information" on $A^\circ$. As for the morphisms between such objects, we should consider coalgebra homomorphisms $f:A^\circ\rightarrow B^\circ$ such that $f^{-1}(V)\in Q_A$ for all $V\in Q_B$. 

\begin{de}\label{rc}
A triple $(C, Q, \Sec)$ will be called a \emph{ringed coalgebra} if $C$ is a coalgebra, $Q$ is a subset of $P_C$ such that 
\begin{itemize} 
\item $Q$ is closed under finite wedge product,
\item $Q$ is closed under (not necessarily finite) intersection,
\item $Q$ contains $C$ and $0$, \end{itemize}
and if $\Sec$ is a subfunctor of the dual functor $(-)^*:P_C^{op}\rightarrow\rfd$, i.e., $\Sec(D)\subset D^*$ for every subcoalgebra $D\subset C$, and the following diagram commutes:
\[\begin{tikzcd}
	{D_2^*} & {D_1^*} \\
	{\Sec(D_2)} & {\Sec(D_1)}
	\arrow["{i^*}", from=1-1, to=1-2]
	\arrow[hook, from=2-1, to=1-1]
	\arrow["{\rho^{D_2}_{D_1}}"', from=2-1, to=2-2]
	\arrow[hook, from=2-2, to=1-2]
\end{tikzcd}\]
whenever $D_1\subset D_2$. Here $i:D_1\hookrightarrow D_2$ is the inclusion and the algebra homomorphism $\rho^{D_2}_{D_1}:\Sec(D_2)\rightarrow\Sec(D_1)$ is the one induced by the inclusion $D_1\subset D_2$ and will be called a restriction from $D_2$ to $D_1$.

The coalgebra $C$ and the collection $Q$ of subcoalgebras will respectively be called the \emph{underlying coalgebra} and the \emph{topology} of the ringed coalgebra $(C, Q, \Sec)$. For every subcoalgebra $D\subset C$, the subalgebra $\Sec_C(D)\subset D^*$ will be called the \emph{section} of $D$. We often omit $Q_C$ and $\Sec_C$ when it does not cause any confusion. 

%If $C_1$ and $C_2$ are subcoalgebras of $C$ such that $C_1\subset C_2$, then $\Sec$ gives 
If $C_1=(C_1, Q_1, \Sec_1)$ and $C_2=(C_2, Q_2, \Sec_2)$ are ringed coalgebras, then a \emph{morphism of ringed coalgebras} from $C_1$ to $C_2$ is a coalgebra homomorphism $f:C_1\rightarrow C_2$ satisfying the following:
\begin{itemize}
\item  For every $D\in Q_2$, $f^\dd(D)$ lies in $Q_1$ and the dual map $D^*\rightarrow f^\dd(D)^*$ restricts to an algebra homomorphism $f^*_D:\Sec_2(D)\rightarrow\Sec_1(f^\dd(D))$.
\end{itemize}
We denote by $\cs$ the category of ringed coalgebras and the morphisms.
\end{de}

Every ringed coalgebra $C=(C, Q, \Sec)$ produces an RFD algebra by applying the functor $\Sec:P^{op}_C\rightarrow\rfd$ to the underlying coalgebra $C$. This functor is an analogy of the global section functor on ringed spaces.
\begin{de}
We define a functor $\Gamma:\cs\rightarrow \rfd^{op}$ by
$$\Gamma(C):=\Sec(C)\subset C^*$$
for every ringed coalgebra $C=(C, Q, \Sec)$. This functor $\Gamma$ will be called the global section functor.
\end{de}
Note that there is a natural coalgebra homomorphism $\iota_C:C\rightarrow\Gamma(C)^\circ$ obtained by composing $i_C:C\hookrightarrow C^{*\circ}$ and the canonical surjection $C^{*\circ}\twoheadrightarrow\Gamma(C)^\circ$. 

If $(C, Q, \Sec)$ is a ringed coalgebra, then $pts(C)=\{x\in C|\Delta(x)=x\otimes x, \varepsilon(x)=1\}$ can be endowed with the weakest topology such that a subset of the form $pts(C)\cap D, D\in Q$ is closed. This defines a functor $pts:\cs\rightarrow\tops$, where $\tops$ denotes the category of topological spaces and continuous maps, and it gives the following commutative diagram:
\[\begin{tikzcd}
	\cs & \tops \\
	\cog & \set
	\arrow["pts", from=1-1, to=1-2]
	\arrow[from=1-1, to=2-1]
	\arrow[from=1-2, to=2-2]
	\arrow["pts", from=2-1, to=2-2]
\end{tikzcd}\]
Here the vertical arrows are the forgetful functors.

%\subsection{Topological spaces and $\cs$}

%Our next goal is to extend the adjunction between the categories $\cog$ and $\set$ given by $k[-]$ and $pts$ to the one between the categories $\cs$ and $\tops$.

\section{RFD algebras and ringed coalgebras}
\subsection{Ringed coalgebras associated to fully RFD algebras}\label{cog1}
In this section we construct ringed coalgebras from RFD algebras. For simplicity, we will focus on fully RFD algebras introduced by \cite{R2}. An algebra $A$ is called \emph{fully residually finite-dimensional}(fully RFD) if every finitely generated left $A$-module is a subdirect product of finite dimensional left $A$-modules. In other words, every finitely generated left $A$-module is an RFD $A$-module. The quotients of a fully RFD algebra are RFD and every maximal left ideal is of finite codimension. The paper \cite{R2} that introduced fully RFD algebras shows that every affine Noetherian PI algebras are fully RFD.

%RFD algebras $A$ satisfying the following properties:
%\begin{itemize}
%\item[(F)] Every maximal left ideal of $A$ is of finite codimensional.%, andevery two-sided ideal of $A$ is intersection of finite codimensional ideals.
%\end{itemize}

%Examples of RFD algebras satisfying the above condition include prime affine PI algebras. Other examples are fully RFD algebras discussed in Manny's paper.

We denote by $\frfd$ the full subcategory of $\rfd$ whose objects are fully RFD algebras. %The construction defines a functor from $\rfd^{op}$ to $\cs$. We restrict this functor to a full subcategory of $\rfd^{op}$ to obtain a fully faithful functor. 

Let $A$ be a fully RFD algebra. We set 
$$Q_A:=\{Z^\circ(I)\subset A^\circ|I\subset A:\ \text{a two-sided ideal}\}.$$
This set is closed under arbitrary intersection and under binary wedge product by \ref{AQ}. It is clear that $Q_A$ contains $0$ and $A^\circ$ as its elements.
Therefore $Q_A$ is a topology on the finite dual coalgebra $A^\circ$.

Our next goal is to constructor a subfunctor $\Sec$ of $(-)^*:P_{A^\circ}^{op}\rightarrow\rfd$. 

%We review the Gabriel localization of rings and prime spectra. See \cite{ZS}, \cite{BS}.

\begin{de}\label{sec}
Let $A$ be a fully RFD algebra and $C\subset A^\circ$ be a subcoalgebra. 
For a subcoalgebra $C$, we define a collection $\F_C$ of left ideals of $A$ by
$$\F_C:=\{I\subset A|I\ \text{is a left ideal},\ Z^\circ(I)\cap C=0\}.$$
Furthermore, we define a subset $S_C$ of $C^*$ by
$$S_C:=\{x\in C^*|\exists I\in\F_C\ a_C(I)x\subset a_C(A)\}$$
 where $a_C$ is the composition of the natural algebra homomorphisms $i_A:A\hookrightarrow A^{\circ*}$ and the canonical surjection $A^{\circ*}\twoheadrightarrow C^*$. We denote by $A_C$ the subalgebra of $C^*$ generated by $S_C$:
 $$A_C:=k\langle S_C\rangle\subset C^*.$$
\end{de} 

Note that $S_C$ contains $a_C(A)$ by definition. If $A^\circ\supset C_1\supset C_2$, then
$\F_{C_1}\subset\F_{C_2}$. The dual map $C_1^*\rightarrow C_2^*$ of the inclusion $C_2\hookrightarrow C_1$ induces a monoid homomorphism $S_{C_1}\rightarrow S_{C_2}$ and hence an algebra homomorphism $A_{C_1}\rightarrow A_{C_2}$. Thus the construction of $A_C$ defines a functor $P_{A^\circ}^{op}\rightarrow\rfd$.

The above definition of $A_C$ is inspired by Gabriel localization(see \cite{BS}). If $A$ is a commutative integral domain, we can show that the Gabriel localization of $A$ with the Gabriel filter 
$$\F_U=\{I\subset A\mid I:\text{ideal},\ Z(I)\cap U=\emptyset\}$$
defined by an open set $U\subset\spc(A)$
is isomorphic to the section $\O_X(U)$ on $U$ where $X=\spc(A)$. We will show in \ref{dsec} that if $A$ is a finitely generated algebra,
there is an isomorphism
$A_{\dist(U)}\simeq\O_X(U)$ where $\dist(U)\subset A^\circ$ is the distribution coalgebra of $U$ explained in \ref{ucos}.

We prepare the following lemma to compute some examples of $A_C$.

\begin{lem}\label{frfd}
Let $A$ be a fully RFD algebra.
\begin{enumerate}
\item[(1)] Every two-sided ideal of $A$ is an intersection of two-sided ideals of finite codimension.
\item[(2)] If a left ideal $I\subset A$ satisfies $Z^\circ(I)=0$, then $I=A$. 
\end{enumerate}
\end{lem}
\begin{proof}
(1) Let $I\subset A$ be a two-sided ideal. Then $A/I$ be a finitely generated left $A$-module. By definition, this is a subdirect product of finite dimensional left modules. This means that $I$ is an intersection of left ideals of finite codimension. The largest two-sided ideal contained in a left ideal of finite codimension has finite codimension. Therefore $I$ is an intersection of two-sided idels of finite codimension.\\
(2) We show the contraposition. Take a left ideal $I\subsetneq A$. There exists a maximal left ideal $J\subsetneq A$ containing $I$. By assumption $J$ has finite codimension and hence 
so does the algebra of endomorphisms on the left $A$-module $A/J$. The kernel of the algebra homomorphism $A\rightarrow End(A/J)$ is the two-sided ideal 
$$\ann(A/J)=\{x\in A|xA\subset J\}\subset J$$
so that $\ann(A/J)$ is a two-sided ideal of finite codimension contained in $J$.
Thus any $\phi\in A^*$ whose kernel contains $J$ is in $Z^\circ(J)$. Therefore $Z^\circ(J)\neq0$ and hence $Z^\circ(I)\neq0$.
\end{proof}

\exam\label{ac} Let $C=Z^\circ(I)$ for some two-sided ideal $I\subset A$ of an RFD algebra $A$ in $\frfd$. Then we have $A_C=S_C=a_C(A)\simeq A/I$. Indeed, let $\phi\in S_C$. By the definition of $S_C$, $a_C(J)\phi\subset a_C(A)$ for some left ideal $J\subset A$ such that $Z^\circ(J)\cap Z^\circ(I)=0$. Here $Z^\circ(J)\cap Z^\circ(I)=Z^\circ(I+J)$ so that $Z^\circ(I+J)=0$. This is equivalent to $J+I=A$ by \ref{frfd}. Since 
$$a_C(J)\phi=a_C(I+J)\phi = a_C(A)\phi,$$
$\phi\in S_C$ implies $\phi\in a_C(A)$. This shows $S_C=a_C(A)$.
Recall that $a_C$ is the composite of the inclusion $A\hookrightarrow A^{\circ*}$ and the canonical surjective map $A^{\circ*}\rightarrow C^*$. Therefore an element $a\in A$ is in $\ker a_C$ if and only if $\phi(a)=0$ for all $\phi\in C=Z^\circ(I)$, i.e., $a\in Z^\circ(I)^\perp$. By \ref{frfd}, $I$ is an intersection of two-sided ideals of finite codimension. Therefore it is closed subspace of $A$ in the sense of \ref{os} and hence $Z^\circ(I)^\perp=I$. Thus $a_C(A)\simeq A/I$.

\exam\label{ac2} If every ideal in $\F_C$ is two-sided(e.g. $A$ is commutative), then $S_C=A_C$. Indeed, if $x, y\in C^*$ satisfies $a_C(I)x, a_C(J)y\subset a_C(A)$ for some $I, J\in \F_C$, then
$$a_C(JI)xy\subset a_C(JA)y=a_C(J)y\subset a_C(A)$$
and
$$a_C(IJ)(x+y)\subset a_C(I\cap J)(x+y)\subset a_C(I\cap J)x+a_C(I\cap J)y\subset a_C(A).$$
Here the product $IJ$ is in $\F_C$ since
$$Z^\circ(IJ)\cap C=(Z^\circ(I)\vee Z^\circ(J))\cap C=(Z^\circ(I)\cap C)\vee_C (Z^\circ(J)\cap C)=0\vee_C 0=0.$$
Therefore $xy, x+y\in S_C$. The subset $S_C\subset C^*$ contains the additive identity $0$ and the multiplicative identity $\eps_C$ so it is a subalgebra of $C^*$. Therefore the subset $S_C\subset C^*$ is a subalgebra so that $A_C=S_C$.

We will see more concrete examples for commutative case(see \ref{dsec}).

\begin{de}\label{algrc}
Let $A$ be a fully RFD algebra. We define a ringed coalgebra associated to $A$ to be the triple $(A^\circ, Q_A, \Sec_A)$ where 
$$Q_A:=\{Z^\circ(I)\subset A^\circ|I\subset A:\ \text{a two-sided ideal}\}$$
and 
$$\Sec_A(C):=A_C\ \text{for every subcoalgebra}\ C\subset A^\circ.$$
Abusing notation, we often denote by $A^\circ$ the triple $(A^\circ, Q_A, \Sec_A)$ if it does not cause any confusion.
\end{de}

We show the functoriality of the construction. Recall that if $\varphi:A\rightarrow B$ is an algebra homomorphism and $f=\varphi^\circ:B^\circ\rightarrow A^\circ$ is the finite dual of $\varphi$, then by \ref{algpb} we have
$$f^\dd(Z^\circ(I))=Z^\circ(B\varphi(I)B)$$
for all ideal $I\subset A$. In particular, for every $C\in Q_A$, $f^\dd(C)\in Q_B$.

\begin{lem}
Let $\varphi:A\rightarrow B$ be an algebra homomorphism, $f:=\varphi^\circ$ be its finite dual, and $C\subset A^\circ$ be a subcoalgebra. 
\begin{enumerate}
\item The assignment $I\mapsto B\varphi(I)$ defines a map from $\F_C$ to $\F_{f^\dd(C)}$.
\item The dual of the restriction $f|_{f^\dd(C)}:f^\dd(C)\rightarrow C$ of $f$ restricts to a map from $S_C$ to $S_{f^\dd(C)}$ and hence an algebra homomorphism from $A_C$ to $B_{f^\dd(C)}$.
\end{enumerate}
\end{lem}
\begin{proof}
(1) Let $I\in \F_C$. It suffices to show that $Z^\circ(B\varphi(I))\cap f^\dd(C)=0$. Let $\phi\in Z^\circ(B\varphi(I))\cap f^\dd(C)$. Then $f(\phi)=\phi\circ\varphi\in Z^\circ(I)\cap C$ so that $f(\phi)=0$. This shows that $ Z^\circ(B\varphi(I))\cap f^\dd(C)\subset\ker f$. In particular $ Z^\circ(B\varphi(I))\cap f^\dd(C)\subset\ker\eps$ since $\ker f\subset\ker\eps$ where $\eps$ is the counit of $B^\circ$. Since $Z^\circ(B\varphi(I))\cap f^\dd(C)$ is a left coideal, we have $\Delta(\phi)=\sum\psi_{i, 1}\otimes \psi_{i, 2}$ for some finite elements $\psi_{i, 1}\in  Z^\circ(B\varphi(I))\cap f^\dd(C)$ and $\psi_{i, 2}\in B^\circ$ where $1\leq i\leq n$. Therefore
$$\phi=(\eps\otimes id)(\Delta(\phi))=\sum_{1\leq i\leq n}\eps(\psi_{i, 1})\psi_{i, 2}=\sum_{1\leq i\leq n}0\cdot\psi_{i, 2}=0.$$
(2) Let $g$ be the dual algebra map $C^*\rightarrow f^\dd(C)^*$ of $f|_{f^\dd(C)}$. We show that the dual map $g$ restricts to a map from $S_C$ to $S_{f^{\dd}(C)}$. Note that there is a unique algebra homomorphism $\varphi':a_C(A)\rightarrow a_{f^{\dd}(C)}(B)$ that makes the following diagram commute: 
\[\begin{tikzcd}
	{A^{\circ*}} && {B^{\circ*}} \\
	& {C^*} && {f^\dd(C)^*} \\
	A && B \\
	& {a_C(A)} && {a_{f^\dd(C)}(B)}
	\arrow["{\varphi^{\circ*}}", from=1-1, to=1-3]
	\arrow[two heads, from=1-1, to=2-2]
	\arrow[two heads, from=1-3, to=2-4]
	\arrow["g"{pos=0.3}, from=2-2, to=2-4]
	\arrow["i_A", hook, from=3-1, to=1-1]
	\arrow["\varphi"{pos=0.3}, from=3-1, to=3-3]
	\arrow[two heads, from=3-1, to=4-2]
	\arrow["i_B"{pos=0.3}, hook, from=3-3, to=1-3]
	\arrow[two heads, from=3-3, to=4-4]
	\arrow[hook, from=4-2, to=2-2]
	\arrow["\varphi'", dashed, from=4-2, to=4-4]
	\arrow[hook, from=4-4, to=2-4]
\end{tikzcd}\]

Let $x\in S_C$ so that $a_C(I)x\subset a_C(A)$ for some $I\in \F_C$. Then 
$$a_{f^\dd(C)}(\varphi(I))g(x)=g(a_C(I)x)\subset g(a_C(A))\subset a_{f^\dd(C)}(B)$$ 
and hence $a_{f^\dd(C)}(B\varphi(I))g(x)\subset a_{f^\dd(C)}(B)$. Since $B\varphi(I)\in\F_{f^\dd(C)}$, $g(x)\in S_{f^\dd(C)}$. The proof is complete.
\end{proof}

%In the usual algebraic geometry, the construction of affine varieties define a functor and the global sections of an affine variety $X=\spm(A)$ recovers the algebra $A$. The following is an analogue in our setting:

\begin{de}
The assignment that sends a fully RFD algebra $A$ to the ringed coalgebra $A^\circ=(A^\circ, Q_A, \Sec_A)$ defines a functor $\frfd^{op}\rightarrow\cs$. By abusing notation, we denote the functor also by $(-)^\circ$.
\end{de}

We have associated to each fully RFD algebra $A$ a ringed coalgebra $A^\circ$. On the other hand, each ringed coalgebra $(C, Q, \Sec)$ gives an RFD algebra by the global section functor $\Gamma:\cs\rightarrow\rfd^{op}$ defined by $\Gamma(C)=\Sec(C)$. This functor recovers the original algebra $A$ from the associated ringed coalgebra $A^\circ$.

\begin{cor}
 Let $A$ be a fully RFD algebra. Then $\Gamma\circ(-)^\circ:\frfd^{op}\rightarrow\frfd^{op}$ is equal to the identity functor $\id_{\frfd^{op}}$ on $\frfd^{op}$.
\end{cor}
\begin{proof}
The equality $\Gamma(A^\circ)=A$ follows from \ref{ac}:
$$\Gamma(A^\circ)=A_{A^\circ}=a_{A^\circ}(A)=A.$$
The equality $\Gamma\circ (-)^\circ=\id_{\frfd^{op}}$ on the morphisms of $\frfd^{op}$ follows from the fact that the diagram
\[\begin{tikzcd}
	{A^{\circ*}} & {B^{\circ*}} \\
	A & B
	\arrow["{f^{\circ*}}", from=1-1, to=1-2]
	\arrow[hook, from=2-1, to=1-1]
	\arrow["f", from=2-1, to=2-2]
	\arrow[hook, from=2-2, to=1-2]
\end{tikzcd}\]
commutes for all algebra homomorphism $f:A\rightarrow B$.
\end{proof}

%\begin{cor}
%Let $A$ be a fully RFD algebra and $A^\circ$ be the associated ringed coalgebra. Then $\Gamma(A^\circ)=A$.
%\end{cor}
%\begin{proof}
%We have seen in \ref{ac} that $A_{Z^\circ(I)}=A/I$ for any two-sided ideal $I\subset A$. Since $A^\circ=Z^\circ((0))$, we have
%$$\Gamma(A^\circ)=A_{A^\circ}=A_{Z^\circ((0))}=A/(0)=A.$$
%\end{proof}

Next, we define a full subcategory $\csu$ of $\cs$ whose objects $(C, Q, \Sec)$ satisfy the following property:
$$ \text{(A)\ The homomorphism $\iota_C$ is a morphism from $(C, \Sec_C)$ to $(\Gamma(C)^\circ, \Sec_{\Gamma(C)})$ in $\cs$.}$$
Note that the ringed coalgebra $A^\circ=(A^\circ, \Sec_A)$ for some algebra $A$ in $\frfd$ satisfies (A). Indeed, we have seen that $\Gamma(A^\circ)=A$. The canonical coalgebra homomorphism 
$\iota_{A^\circ}:A^\circ\rightarrow\Gamma(A^\circ)^\circ=A^\circ$ is nothing but the identity morphism of $A^\circ$.

\begin{thm}\label{algadj}
Let $A$ be a fully RFD algebra and $C$ be a ringed coalgebra satisfying (A). There is a natural bijective correspondence
$$\csu(C, A^\circ)\simeq \rfd(A, \Gamma(C))$$
given by  $f\mapsto \Gamma(f)$. 
\end{thm}
\begin{proof}
We show that the assignment $g\mapsto g^\circ\circ\iota_C$ gives the inverse of $f\mapsto \Gamma(f)$. The following commutative diagram shows that $f=\Gamma(f)^\circ\circ\iota_C$:
\[\begin{tikzcd}
	C & {C^{*\circ}} & {\Gamma(C)^\circ} \\
	{A^\circ} & {A^{\circ*\circ}} & {\Gamma(A^\circ)^\circ=A^\circ}
	\arrow["i_C", hook, from=1-1, to=1-2]
	\arrow["f"', from=1-1, to=2-1]
	\arrow[from=1-2, to=1-3]
	\arrow["{f^{*\circ}}"', from=1-2, to=2-2]
	\arrow["{\Gamma(f)^\circ}"', from=1-3, to=2-3]
	\arrow["i_{A^\circ}", hook, from=2-1, to=2-2]
	\arrow[from=2-2, to=2-3]
\end{tikzcd}\]
Here the bottom morphism is identity on $A^\circ$. We consider the following commutative diagram in order to show $g=\Gamma(g^\circ\circ\iota_C)$:
\[\begin{tikzcd}
	& {A^{\circ*}} \\
	{\Gamma(A^\circ)} & A & {\Gamma(C)^{\circ*}} \\
	& {\Gamma(\Gamma(C)^\circ)} & {\Gamma(C)}
	\arrow["{g^{\circ*}}"{pos=0.2}, from=1-2, to=2-3]
	\arrow[hook, from=2-1, to=1-2]
	\arrow["\id_A"', from=2-1, to=2-2]
	\arrow["{\Gamma(g^\circ)}"', from=2-1, to=3-2]
	\arrow["i_A"', hook, from=2-2, to=1-2]
	\arrow["g"{pos=0.1}, from=2-2, to=3-3]
	\arrow[hook, from=3-2, to=2-3]
	\arrow["{\Gamma(\iota_C)}"', from=3-2, to=3-3]
	\arrow["i_{\Gamma(C)}"', hook, from=3-3, to=2-3]
\end{tikzcd}\]
The rightmost square and the tilted square are commutative. The injectivity of the natrual map $\Gamma(C)\rightarrow\Gamma(C)^{\circ*}$ implies the commutativity of the bottom square. Since $\Gamma(\iota_{A^\circ})=\id_A$, we have $g=\Gamma(g^\circ\circ\iota_C)$.
\end{proof}

In particular, 
$$\Gamma(C)\simeq\rfd(k[x], \Gamma(C))\simeq\csu(C, k[x]^\circ)$$ 
as sets if $C$ satisfies condition (A). 

\begin{cor}\label{algemb}
The functor $(-)^\circ:\frfd^{op}\rightarrow\csu$ is fully-faithful.
\end{cor}
\begin{proof}
It follows from the natural isomorphisms
$$\csu(B^\circ, A^\circ)\simeq\frfd(A, \Gamma(B^\circ))\simeq\frfd(A, B)$$
for any algebras $A$ and $B$ in $\frfd$.
\end{proof}

If $A$ is commutative and affine, the subset of the $k$-rational points of $\spc(A)$ is the set $\spm(A)$ of algebra homomorphism from $A$ to $k$ and it is endowed with the Zariski topology. The restriction of this functor to the category $\aff$, the full subcategory of $\calg$ whose objects are commutative affine algebras, makes the following diagram commute up to natural isomorphism:
\[\begin{tikzcd}
	\aff^{op} && \csu \\
	\\
	& \tops
	\arrow["{(-)^\circ}", hook, from=1-1, to=1-3]
	\arrow["\spm"', from=1-1, to=3-2]
	\arrow["pts", from=1-3, to=3-2]
\end{tikzcd}\]

\section{Schemes locally of finite type and ringed coalgebras} 

We have associated a ringed coalgebra to each fully RFD algebra by using the finite dual coalgebra. These ringed coalgebras can be thought of as noncommutative generalization of affine schemes over a fixed field $k$. In this section we associate 
ringed coalgebras to schemes locally of finite type over $k$. This is done by endowing the underlying coalgebra $\dist(X)$ of a scheme $X$ with the topology $Q_X$ obtained from the closed subschemes of $X$ and with the  largest algebras $\Sec_X(C)$ of sections on subcoalgebras $C\subset\dist(X)$ where every scheme morphism from an affine scheme to $X$ gives a morphism to $\dist(X)$ of ringed coalgebras. The construction defines a faithful functor that is compatible with the fully-faithful functor $(-)^\circ:\frfd^{op}\rightarrow\csu$ in \ref{algrc} and with the functor $\lv -\rv:\sch^{lf}\rightarrow\tops$ that gives the sets of closed points of the underlying topological spaces. Furthermore, we show that it is also full if the base field $k$ is algebraically closed. 

Throughout this section all algebras are commutative RFD and all schemes are the schemes over $k$. We say a scheme is locally of finite type if it is locally of finite type over $k$. For a scheme $X$ we denote by $\lv X\rv$ the set of points of $X$ whose residue field is a finite extension of $k$. We denote by $\sch$(resp. $\sch^{lf}$) the category of schemes(resp. schemes locally of finite type) over $k$. For the fundamentals, see \cite{H}.

\subsection{Ringed coalgebras associated to schemes locally of finite type}\label{cog2}
%In this subsections, all algebras are commutative and all schemes are locally of finite type over $k$. Our goal in this subsection is to define a ringed coalgebra $\emb(X)$ associated to a scheme $X$ locally of finite type. %We will assume that $k$ is algebraically closed.

Let $X$ be a scheme locally of finite type. We equip the underlying coalgebra $\dist(X)$ with a collection $Q_X$ of subcoalgebras of $\dist(X)$ and a functor $\Sec_X:P^{op}_{\dist(X)}\rightarrow\rfd$ to assign a ringed coalgebra to $X$. 

Recall that iff $X$ is affine, i.e. $X=\spc(A)$ for some finitely generated algebra $A$, then there is a correspondence between the closed subschemes of $X$ and the ideals of $A$. Every ideal $I\subset A$ defines a subcoalgebra $Z^\circ(I)\subset A^\circ$ and the topology $Q_A$ is defined to be the set of such subcoalgebras. Similarly, if $i:Z\rightarrow X$ is a closed immersion, then the homomorphism at each stalk $\O_{X, i(z)}\rightarrow\O_{Z, z}$ is surjective and hence the finite dual $\O_{Z, z}^\circ\rightarrow\O_{X, i(z)}^\circ$ is injective. Thus the induced coalgebra homomorphism $\dist(i):\dist(Z)\rightarrow\dist(X)$ is injective and $\dist(Z)$ can be seen as a subcoalgebra of $\dist(X)$. The image of $\dist(i)$ only depends on the isomorphism class of closed immersions to $X$, so every closed subscheme of $Z$ defines a subcoalgera $\dist(Z)$. We set
%$$Q_X:=\{Im(\dist(i))\subset\dist(X)|i:Z\hookrightarrow X\ \text{is a closed immersion}\}.$$
$$Q_X:=\{\dist(Z)\subset\dist(X)|Z\ \text{is a closed subscheme of}\ X\}.$$
This collection $Q_X$ is closed under intersection and finite wedge product. If $\{Z_i\hookrightarrow X\}$ is a collection of closed immersions and $\{I_i\}$ are the corresponding ideal sheaves, then the closed subscheme $\bigcap_i Z_i$ corresponding to the ideal sheaf $\sum_iI_i$, which is called the scheme theoretic intersection, gives
$$\dist\left(\bigcap_i Z_i\right)=\bigcap_i \dist(Z_i)$$
since $\dist$ preserves pullbacks. For closed immersions $Z_1, Z_2\hookrightarrow X$ with corresponding ideals sheaves $I$ and $J$ on $X$, consider the closed subscheme $Z_{12}$ corresponding to the ideal sheaf $IJ$. By \ref{sumprod} we have
\begin{eqnarray*}
\dist(Z_{12})&=&\bigoplus_{x\in X}Z^\circ(I_xJ_x)=\bigoplus_{x\in X}\left(Z^\circ(I_x)\vee Z^\circ(J_x)\right)\\
&=&\left(\bigoplus_{x\in X}Z^\circ(I_x)\right)\vee\left(\bigoplus_{x\in X}Z^\circ(J_x)\right)=\dist(Z_1)\vee\dist(Z_2)
\end{eqnarray*}
where $Z^\circ(-)$ stands for the vanishing subspaces of $\O_{X, x}^\circ$. Thus $Q_X$ is closed under wedge products. It is clear that $Q_X$ contains $0$ and $\dist(X)$. Therefore $Q_X$ defines a topology on the coalgebra $\dist(X)$.

Next, we define a functor $\Sec_X:P^{op}_{\dist(X)}\rightarrow \rfd$. We aim to do so in such a way that if $f:\spc(A)\rightarrow X$ is a morphism of scheme , then $g=\dist(f):A^\circ\rightarrow \dist(X)$ is a morphism of ringed coalgebras. Recall that for a subcoalgebra $C\subset\dist(X)$, $g^\dd(C)$ is the largest subcoalgebra of $A^\circ$ that is contained in the preimage $g^{-1}(C)$.
We denote by $g_C:g^\dd(C)\rightarrow C$ the restriction of $g$ to $g^\dd(C)$. The coalgebra homomorphism $g$ is a morphism of ringed coalgebras if the dual $g_C^*:C^*\rightarrow g^\dd(C)^*$  restricts to an algebra homomorphism $\Sec_X(C)\rightarrow A_{g^\dd(C)}$, i.e., there is an algebra homomorphism $\Sec_X(C)\rightarrow A_{g^\dd(C)}$ that makes the following diagram commute:
\[\begin{tikzcd}
	{C^*} & {g^\dd(C)^*} \\
	{\Sec_X(C)} & {\Sec_A(g^\dd(C))=A_{g^\dd(C)}}
	\arrow["g_C^*", from=1-1, to=1-2]
	\arrow["incl.", hook, from=2-1, to=1-1]
	\arrow["\exists", dashed, from=2-1, to=2-2]
	\arrow["incl.", hook, from=2-2, to=1-2]
\end{tikzcd}\]

\begin{de}\label{secdef}
Let $X$ be a scheme locally of finite type and $C\subset\dist(X)$ be a subcoalgebra. We define $\Sec_X(C)\subset C^*$ to be the intersection
$$\Sec_X(C)=\bigcap_{g=\dist(f)}g_C^{*-1}(A_{g^\dd(C)})\subset C^*$$
of the preimages of $A_{g^\dd(C)}$
where $f$ runs over all scheme morphisms $\spc(A)\rightarrow X$ from the spectra of finitely generated algebras $A$ to $X$ and $g_C:g^\dd(C)\rightarrow C$ is the coalgebra homomorphism induced by restricting $g=\dist(f):A^\circ\rightarrow\dist(X)$ to $g^\dd(C)$. 
\end{de}

Equivalently, $\Sec_X(C)\subset C^*$ is the subalgebra generated by elements $\phi\in C^*$ satisfying the following condition:
\begin{itemize}
\item[(S)] For any commutative finitely generated algebra $A$ and a scheme morphism $f:\spc(A)\rightarrow X$, there exists an ideal $I\subset A$ such that $Z^\circ(I)\cap g^\dd(C)=0$ and $a_{g^\dd(C)}(I)g^*(\phi)\subset a_{g^\dd(C)}(A)$ where $g=\dist(f)$.
\end{itemize}
%$$\Sec_X(C)=\{\phi\in C^*|\exists U\subset X:\text{open}\ s.t.\ C\subset\dist(U) \exists \varphi\in\O(U)\ \phi=\varphi|_C)\}.$$

%The above condition means that the dual map $C^*\rightarrow g^\dd(C)^*$ of the coalgebra homomorphism $g^\dd(C)\rightarrow C$ restricts to an algebra homomorphism from $\Sec_X(C)$ to $\Sec_A(g^\dd(C))=A_{g^\dd(C)}$ whenever $g:A^\circ\rightarrow\dist(Y)$ is induced by some scheme morphism $\spc(A)\rightarrow Y$:
%\[\begin{tikzcd}
	%{C^*} & {g^\dd(C)^*} \\
	%{\Sec_X(C)} & {\Sec_A(g^\dd(C))=A_{g^\dd(C)}}
	%\arrow["g_C^*", from=1-1, to=1-2]
	%\arrow["incl.", hook, from=2-1, to=1-1]
	%\arrow[dashed, from=2-1, to=2-2]
	%\arrow["incl.", hook, from=2-2, to=1-2]
%\end{tikzcd}\]
%where $g_C:g^\dd(C)\rightarrow C$ is the coalgebra homomorphism induced by restricting $g:A^\circ\rightarrow \dist(X)$ to $g^\dd(C)$. 

If $\dist(X)\supset C_1\supset C_2$ are subcoalgebras, then the dual map $i^*:C_1^*\rightarrow C_2^*$ restricts to an algebra homomorphism $\Sec_X(C_1)\rightarrow\Sec_X(C_2)$. In fact, let $f:\spc(A)\rightarrow X$ be a morphism of schemes where $A$ is an affine algebra and $g:=\dist(f)$. Then the squares on the top, in front and behind in the following diagram commute:
\[\begin{tikzcd}
	{C_1^*} & {g^\dd(C_1)^*} \\
	{g_C^{*-1}(A_{g^\dd(C_1)})} & {g_C^{*-1}(A_{g^\dd(C_1)})} && {C_2^*} & {g^\dd(C_2)^*} \\
	&&& {A_{g^\dd(C_1)}} & {A_{g^\dd(C_2)}}
	\arrow["{g_{C_1}^*}", from=1-1, to=1-2]
	\arrow["{i^*}"{pos=0.6}, shift right, from=1-1, to=2-4]
	\arrow["{i^*}", from=1-2, to=2-5]
	\arrow["{incl.}", hook, from=2-1, to=1-1]
	\arrow[from=2-1, to=2-2]
	\arrow[dashed, from=2-1, to=3-4]
	\arrow["{incl.}"{pos=0.3}, hook, from=2-2, to=1-2]
	\arrow[from=2-2, to=3-5]
	\arrow["{g_{C_2}^*}"', from=2-4, to=2-5]
	\arrow["{incl.}"', hook, from=3-4, to=2-4]
	\arrow[from=3-4, to=3-5]
	\arrow["{incl.}"', hook, from=3-5, to=2-5]
\end{tikzcd}\]
%the dual of an inclusion of coalgebras restricts to an algebra homomorphism between their sections. 
Here $i$ stands for an inclusion. We have seen in section 3 that $i^*:g^\dd(C_1)^*\rightarrow g^\dd(C_2)^*$ restricts to an algebra homomorphism $\Sec_{g^\dd(C_1)}\rightarrow\Sec_{g^\dd(C_2)}$. Therefore $C_1^*\rightarrow C_2^*$ restricts to 
$g_C^{*-1}(\Sec_{g^\dd(C_1)})\rightarrow g_C^{*-1}(\Sec_{g^\dd(C_2)})$.
Since $f$ is arbitrary, $C_1^*\rightarrow C_2^*$ restricts to an algebra homomorphism
$$\Sec_X(C_1)=\bigcap_{g=\dist(f)}g_{C_1}^{*-1}(A_{g^\dd(C_1)})\rightarrow\bigcap_{g=\dist(f)}g_{C_2}^{*-1}(A_{g^\dd(C_2)})=\Sec_X(C_2).$$
Therefore the definition of $\Sec_X(C)$ for every subcoalgebra $C\subset\dist(X)$ defines a functor $\Sec_X:P^{op}_{\dist(X)}\rightarrow \rfd$.

%In fact, if $Z^\circ(I)\cap g^\dd(C_1)=0$ and $a_{g^\dd(C_1)}(I)g^*(\phi)\subset a_{g^\dd(C_1)}(A)$, \todo{then it also holds that} $Z^\circ(I)\cap g^\dd(C_2)=0$ and $a_{g^\dd(C_2)}(I)g^*(\phi)\subset a_{g^\dd(C_2)}(A)$. 

\begin{de}
We define a ringed coalgebra associated to $X$ to be a triple $(\dist(X), Q_X, \Sec_X)$. Abusing notation, we often denote the triple by $\dist(X)$ if it does not cause any confusion.
\end{de}

Our next goal is to show that the above definition of ringed coalgebras $\dist(X)$ gives a functor from $\sch^{lf}$ to $\cs$. First, we need to check that for a morphism $f:X\rightarrow Y$ of schemes, the induced coalgebra homomorphism $g=\dist(f)$ pulls back every element of $Q_Y$ to an element of $Q_X$. Recall that $\ccog$ denotes the full subcategory of $\cog$ whose objects are cocomutative coalgebras.

\begin{lem}
Let $f:X\rightarrow Y$ be a morphism of schemes and $g:=\dist(f):\dist(X)\rightarrow \dist(Y)$ be the coalgebra homomorphism associated to $f$. Let $D=\dist(Z)\subset \dist(Y)$ be a subcoalgebra defined by a closed subscheme $Z\subset Y$. Then $g^\dd(D)=\dist(f^{-1}(Z))$ where $f^{-1}(Z)\subset X$ is the basechange of $Z$ along $f$. 
\end{lem}
\begin{proof}
By definition the functor $\dist:\sch\rightarrow\ccog$  preserves limits. Thus the subcoalgebra $\dist(f^{-1}(Z))\subset \dist(X)$ is a category theoretic pull-back of the following square on the right in $\ccog$:
\[\begin{tikzcd}
	X && Y && {\dist(X)} && \dist(Y) \\
	\\
	{f^{-1}(Z)=X\times_YZ} && Z && {g^\dd(D)=\dist(X)\times_{\dist(Y)}\dist(Z)} && \dist(Z)
	\arrow["f", from=1-1, to=1-3]
	\arrow["g", from=1-5, to=1-7]
	\arrow["closed", from=3-1, to=1-1]
	\arrow[from=3-1, to=3-3]
	\arrow["closed", from=3-3, to=1-3]
	\arrow[hook, from=3-5, to=1-5]
	\arrow[from=3-5, to=3-7]
	\arrow[hook, from=3-7, to=1-7]
\end{tikzcd}\]
On the other hand, the subcoalgebra $g^\dd(D)\subset\dist(X)$ is a a category theoretic pull-back of the square on the right in $\cog$ by \ref{pb} and hence in $\ccog$. Thus we must have $g^\dd(D)=\dist(f^{-1}(Z))$.
\end{proof}

\begin{lem}\label{tfunc}
Let $f:X\rightarrow Y$ be a morphism of schemes. Then $\alpha:=\dist(f):\dist(X)\rightarrow\dist(Y)$ is a morphism of ringed coalgebras.
\end{lem}
\begin{proof}
Let $C\subset\dist(Y)$ be a subcoalgebra. We need to show that the dual map $C^*\rightarrow \alpha^\dd(C)^*$ restricts to an algebra homomorphism $\Sec_Y(C)\rightarrow\Sec_X(\alpha^\dd(C))$. Let $\phi\in\Sec_Y(C)$ and let $g:\spc(A)\rightarrow X$ be a morphism of schemes. Set $\beta:=\dist(g)$ and $\gamma:=\dist(f\circ g)$. Since $f\circ g:\spc(A)\rightarrow Y$ is a morphism of scheme, there exists $I\subset A$ such that $Z^\circ(I)\cap \gamma^\dd(C)=0$ and $a_{\gamma^\dd(C)}(I)\gamma^*(\phi)\subset a_{\gamma^\dd(C)}(A)$. Here $\gamma^\dd(C)=\beta^\dd(\alpha^\dd(C))$ 
and $\gamma^*=\beta^*\circ\alpha^*$. Therefore $\alpha^*(\phi)$ satisfies condition (S) and is an element of $\Sec_X(\alpha^\dd(C))$.
\end{proof}

\begin{cor}\label{tfaith1}
The assignment that sends a $k$-scheme $X$ locally of finite type to the ringed algebra $\dist(X)=(\dist(X), Q_X, \Sec_X)$ defines a functor $\emb:\sch^{lf}\rightarrow\cs$. Furthermore, $\emb$ is faithful.
\end{cor}
\begin{proof}
The first assertion follows from the preceding lemma. The second assertion follow from the faithfulness of the functor $\dist:\sch^{lf}\rightarrow\cog$ shown in \cite{T2}.
\end{proof}

We can prove the faithfulness of $\emb$ also by using the quasi-separation of schemes locally of finite type. Note that for any two points in a quasi-separated scheme, there exists an affine open subset that contains the two points(\cite[\href{https://stacks.math.columbia.edu/tag/01ZU}{Tag 01ZU}]{SAO}). 

\begin{prop}\label{tfaith2}
The functor $\emb$ is faithful, i.e., the map
$$\emb:\sch^{lf}(X, Y)\rightarrow\cs(\emb(X), \emb(Y))$$
is injective for all schemes $X$ and $Y$ locally of finite type. 
\end{prop}
\begin{proof}
Let $f, g:X\rightarrow Y$ be morphisms of $\sch^{lf}$ such that $\emb(f)=\emb(g)$. Since the subset of closed points of $X$ is dense, it suffices to show that for every closed point $x\in X$, $f$ and $g$ agree on some affine open neighborhood of $x$. Let $x\in X$. Since $Y$ is quasi-separated, there exists an affine open subset $U\subset Y$ that contains $f(x)$ and $g(x)$. Then $f^{-1}(U)\cap g^{-1}(U)$ is an open neighborhood of $x$. Take an affine open neighborhood $V$ of $x$ that is contained in $f^{-1}(U)\cap g^{-1}(U)$. The the restrictions $f|_V, g|_V:V\rightarrow U$ induce $\emb(f|_V), \emb(g|_V):\emb(V)\rightarrow\emb(U)$. However, by assumption we have
$$\emb(f|_V)=\emb(f\circ\rho^X_V)=\emb(f)\circ(\rho^{X}_{V})^\circ=\emb(f)|_{\dist(V)}$$
where $\rho^X_V:\O_X(X)\rightarrow\O_X(V)$ is the restriction of the sections. 
%and$$\rho^{\dist(X)}_{\dist(V)}:\dist(X)\rightarrow\dist(V)$$ stands for a restriction. 
By a similar argument it follows that $$\emb(g|_V)=\emb(g)|_{\dist(V)}.$$
By assumption $\emb(f)|_{\dist(V)}=\emb(g)|_{\dist(V)}$ so that $\emb(f|_V)=\emb(g|_V)$. Moreover $U$ and $V$ are affine so that $\emb(V)=\Gamma(V)^\circ$ and $\emb(U)=\Gamma(U)^\circ$. Thus by the correspondence
$$\csu(\Gamma(V)^\circ, \Gamma(U)^\circ)\simeq\rfd(\Gamma(U), \Gamma(V))\simeq\sch^{lf}(V, U),$$
$\emb(f|_V)=\emb(g|_V)$ implies that $f|_V=g|_V$. The proof is complete.
\end{proof}

We close this section by showing that $\emb$ is compatible with the functors we have introduced so far.

\begin{lem}
Let $A$ be a finitely generated algebra and $X=\spc(A)$. For every subcoalgebra $C\subset A^\circ=\dist(X)$, 
$$\Sec_X(C)=A_C.$$ 
In particular, $\emb(\spc(A))$ and $A^\circ$ are isomorphic as ringed coalgebras. 
\end{lem}
\begin{proof}
Applying condition (S) to the identity morphism on $X=\spc(A)$, we have $\Sec_X(C)\subset A_C$. On the other hand, every element in $A_C$ satisfies condition (S) so that $A_C\subset\Sec_X(C)$.
\end{proof}

\begin{cor}\label{com}
The following diagram commutes up to natural isomorphism:
\[\begin{tikzcd}
	{\aff^{op}} && {\frfd^{op}} \\
	\sch^{lf} && \cs \\
	& \tops
	\arrow["incl.", hook, from=1-1, to=1-3]
	\arrow["\spc"', hook, from=1-1, to=2-1]
	\arrow["{(-)^\circ}", hook, from=1-3, to=2-3]
	\arrow["\emb", from=2-1, to=2-3]
	\arrow["{(-)(k)}"', from=2-1, to=3-2]
	\arrow["pts", from=2-3, to=3-2]
\end{tikzcd}\]
Here $(-)(k):\sch^{lf}\rightarrow\tops$ is the functor that gives the topological space of $k$-rational points with the Zariski topology.
\end{cor}
\begin{proof}
It suffices to show the triangle at the bottom is commutative. For a scheme $X$ locally of finite type, the set of the rational points $X(k)$ and $pts(\emb(X))$ are canonically isomorphic as sets by \cite[2.18]{R2} or \cite[2.1.8]{T2}. The bijection is indeed a homeomorphism since $Z(k)=pts(\emb(Z))$ for every closed subscheme $Z\subset X$ under the identification $X(k)=pts(\emb(X))$.
\end{proof}

If we denote by $\mathbb{A}^1$ the spectrum $\spc(k[x])$, then we have
$$\Gamma(C)\simeq\csu(C, k[x]^\circ)\simeq\csu(C, \emb(\mathbb{A}^1))$$ 
as sets if $C$ satisfies condition (A). The ringed coalgebra $\emb(\mathbb{A}^1)=k[x]^\circ$ can be seen as a affine line in our setting.

\subsection{Schemes locally of finite type over an algebraically closed field}\label{cog3}
We have defined the functor $\emb:\sch^{lf}\rightarrow\cs$ and seen the faithfulness in the previous subsection. In this subsection, \textbf{we assume $k$ to be an algebraically closed field} and show that $\emb$ is full. 

%if it is restricted to the full subcategory $\rsp$ of integral schemes locally of finite type over $k$. 

Let $X$ be a scheme locally of finite type over the algebraically closed field $k$. Therefore the set
$$\lv X\rv=\{x\in X\mid [\kappa(x):k]<\infty\}$$
is equal to the set $X(k)$ of $k$-rational points and this is also the set of closed points in $X$. The subset $\lv X\rv\subset X$ is known to be dense(\cite[\href{https://stacks.math.columbia.edu/tag/01P1}{Tag 01P1}]{SAO}). For every open subset $U\subset X$, there is a natural injective homomorphism
$$\O(U)\hookrightarrow\prod_{x\in U}\O_{X, x}$$
given by sending $f\in\O(U)$ to $(f_x)_{x\in U}\in\prod_{x\in U}\O_{X, x}$ where $f_x$ is the image of $f\in\O(U)$ under the canonical homomorphism $\O(U)\rightarrow\O_{X, x}$. Identfiying $\O(U)$ with a subalgebra of $\prod_{x\in U}\O_{X, x}$, we often write a section $f$ on $U$ as a collection $(f_x)_{x\in U}$ of germs of each point in $U$.

The underlying coalgebra of $X$ can be described as follows(see \ref{ucos}):
$$\dist(X)=\bigoplus_{x\in\lv X\rv}\O_{X, x}^\circ.$$
By \cite[1.1.2(iii)]{T2}, we have
$$\O_{X, x}^{\circ*}=(\dirlim_{n}(\O_{X, x}/\M^n_x)^*)^*\simeq\invlim_n(\O_{X, x}/\M^n_x)=\hat{\O}_{X, x}$$
where $\hat{\O}_{X, x}$ is the completion of the local ring $\O_{X, x}$ by the unique maximal ideal $\M_x$. Therefore
$$\dist(X)^*=\prod_{x\in\lv X\rv}\O_{X, x}^{\circ*}=\prod_{x\in\lv X\rv}\hat{\O}_{X, x}$$
We view the dual algebra of a coalgebra as the algebra of “$k$-valued functions on the coalgebra". The algebra $\O(U)$
of the sections on $U$ can be viewed as a subaglebra of $\dist(U)^*$.
Indeed, we have the following general fact:

\begin{lem}\label{inj}
%Let $X$ be an integral scheme locally of finite type. The algebra homomorphism $\Gamma(X)\hookrightarrow\dist(X)^*$ defined by
%$(x_\p)_\p\mapsto (x_\m)_\m$ is injective.
Let $X$ be a $k$-scheme locally of finite type and $U\subset X$ be an open subset. The composition of the natural homomorphisms
$$\mathcal{O}_X(U)\hookrightarrow\prod_{x\in U}\mathcal{O}_{X, x}\twoheadrightarrow \prod_{x\in \lv U\rv}\mathcal{O}_{X, x}$$
is injective. 
\end{lem}
\begin{proof}
Let $f\in \O_X(X)$ be a global section in the kernel of this homomorphism and $X=\bigcup_i U_i$ be an affine open covering. Then $\lv X\rv=\bigcup_i \lv U_i\rv$. We have the following commutative diagram:
\[\begin{tikzcd}
	{\mathcal{O}_X(X)} && {\prod_{x\in \lv X\rv}\mathcal{O}_{X, x}} \\
	\\
	{\mathcal{O}_X(U_i)} && {\prod_{x\in \lv U_i\rv}\mathcal{O}_{X, x}}
	\arrow[from=3-1, to=3-3]
	\arrow["{rest.}", from=1-1, to=3-1]
	\arrow[two heads, from=1-3, to=3-3]
	\arrow[from=1-1, to=1-3]
\end{tikzcd}\]
Here the horizontal arrow on the bottom is injective. In fact, the arrow can be described as the canonical homomorphism
$$R\rightarrow\prod_{\m\in\spm(R)}R_\m$$
where $R=\O(U_i)$. To show that the homomorphism is injective, let $x$ be an element of the kernel of the homomorphism. Then the annihilator $\ann_R(x)$ cannot be contained in any maximal ideal of $R$. Thus
$1\in\ann_R(x)$ and hence $x=0$. 

By the injectivity of the arrow on the bottom of the diagram, $f|_{U_i}=0$ as an element of $\O_X(U_i)$. By the gluing property of sheaves, we obtain $f=0$ as an element of $\O_X(X)$. 
\end{proof}

Therefore we have a composite of inclusions
$$\O_X(U)\hookrightarrow\prod_{x\in \lv U\rv}\mathcal{O}_{X, x}\hookrightarrow\prod_{x\in \lv U\rv}\hat{\mathcal{O}}_{X, x}=\dist(U)^{*}$$
given by sending $f\in\O_X(U)$ to a collection $\{f_x\}_{x\in \lv U\rv}$.
In this way, the algebra $\O_X(U)$ can be seen as a subalgebra of $\dist(U)^{*}$. Also, the preceding lemma allows us to write every element $f\in \O_X(U)$ as a collection $f=(f_x)_{x\in\lv U\rv}$ of elements $f_x\in\O_{X, x}$. 

%If $X$ is affine, then the image can be characterized as follows:
%\begin{lem}
%Let $X=\spc(A)$ for some finitely generated algebra $A$ and let $U\subset X$ be an open subset. An element 
%$$(f_\m)_{\m\in\lv U\rv}\in\prod_{\m\in\lv U\rv} A_\m$$
%lies in the image of $\O_{X}(U)$ under the injective morphism if and only if for every $\m\in\lv U\rv$ there exists an open neighborhood $V\subset X$ of $\m$ and elements $a, f\in A$ such that for every $\n\in \lv V\rv$ satisfying $f\notin \n$, we have $f_\n=\frac{a}{f}$.
%\end{lem}
%\begin{proof}
%It is immediate from the definition of structure sheaf on $X$ that if $(f_\m)_{\m\in\lv U\rv}$ lies in the image of $\O_X(U)$ then the element satisfies the property. Assume that the element $(f_\m)_{\m\in\lv U\rv}$ satisfies the condition. For $\p\in U$, take an open
%\end{proof}

The following shows that if $A$ is a finitely generated algebra and $$C=\dist(U)\subset A^\circ=\dist(\spc(A))$$ 
is the subcoalgebra associated to $U$, then the subalgebra $\Sec_A(C)\subset C^*$ of sections on a subcoalgebra $C\subset A^\circ$ recovers the section $\O_X(U)$ on $U$.
%The following gives some examples of $A_C$ when $A$ is an affine domain and $C=\dist(U)$ is a subcoalgebra defined by an open set $U\subset\spc(A)$.

\begin{lem}\label{dsec}
Let $A$ be a (commutative) finitely generated algebra and $X=\spc(A)$. Then
$$A_{\dist(U)}\simeq\O_X(U)$$
for every open subscheme $U\subset X$.
\end{lem}
\begin{proof}
We may assume that $U$ is nonempty. By definition, we have
$$A_{\dist(U)}=\{x\in \dist(U)^*|\exists I\in\F_{\dist(U)}\ a_{\dist(U)}(I)x\subset a_{\dist(U)}(A)\}$$
where $a_{\dist(U)}$ is the composite of the following natural algebra homomorphisms:
$$A\hookrightarrow A^{\circ*}\twoheadrightarrow \dist(U)^*=\prod_{\m\in\lv U\rv}\hat{A}_\m.$$

Let $x=(x_\m)_{\m\in \lv U\rv}\in\O_X(U)$ be an element. Since $\spc(A)$ is Noetherian, the open subset $U\subset\spc(A)$ is compact. Take an open covering $U=\bigcup_{1\leq i\leq k} D(f_i)$ so that $x_\p$ can be written as $\frac{a_i}{f^e_i}\in A_{f_i}$ for all $\m\in \lv D(f_i)\rv$. Let $I:=(f_1^{e}, \cdots, f_k^e)$. Note that $I\in \F_{\dist(U)}$. For given $1\leq j\leq k$, we have
$$f^e_j\frac{a_i}{f^e_i}=\frac{a_j}{1}$$
for all $i$ so that $a_{\dist(U)}(I)x\subset a_{\dist(U)}(A)$. Therefore we obtain an injective homomorphism $\O(U)\hookrightarrow A_{\dist(U)}$. On the other hand, suppose that $x\in\dist(U)^*$ satisfies that $a_{\dist(U)}(I)x\subset a_{\dist(U)}(A)$ for some $I\in\F_{\dist(U)}$. We may assume $x\neq0$. Note that $a_{\dist(U)}$ is injective since $U\subset X$ is dense. Therefore we may identify $a_{\dist(U)}(A)$ with $A$. For every $\m\in \lv U\rv$ and an element $f_\m\in I\backslash\m$, we have $f_\m x=a_\m\in A$.  Thus we have
$$f_\n a_\m=f_\n f_\m x=f_\m a_\n$$
for $\m, \n\in\lv U\rv$. We define an element 
$$(x_\p)_\p\in\prod_{\p\in U} A_\p$$
by setting $x_\p=\frac{a_\m}{f_\m}$ if $\p=\m\in\lv U\rv$. If $\p$ is not maximal, then take a maximal ideal $\m\in \lv U\rv$ that contains $\p$ and set $x_\p:=\frac{a_\m}{f_\m}$. This does not depend on the choice of $\m$; indeed, if $\m, \n\in U$ are maximal ideals of $A$ containing $\p\in U$, then $f_\m$ and $f_\n$ are not in $\p$. Since $f_\n a_\m=f_\m a_\n$, we have  
$$\frac{a_\n}{f_\n}=\frac{a_\m}{f_\m}$$
as elements in $A_\p$. We can check that for every $\p\in U$ and $\m\in\lv U\rv$ such that $\p\subset\m$, we have
$$x_\q=\frac{a_\m}{f_\m}$$ 
for all $\q\in U\cap D(f_\m)$.
Therefore the element $(x_\p)_\p\in\prod_{\p\in U} A_\p$ 
is an element of $\O_X(U)$. It is easy to see that the element $(x_\p)_\p$ is mapped to $x\in A_{\dist(U)}$ under the natural injective homomorphism $\O_X(U)\hookrightarrow\dist(U)^*$.
\end{proof}

%Since $\hat{A}_x$ is an integral domain for each $x\in\lv U\rv$, we have $f_\p a_\q-f_\q a_\p=0$ for any $\p, \q\in U$. 
By using the preceding lemma about open subschemes of affine schemes locally of finite type, we can show the following:

\begin{cor}\label{ssec}
Let $X$ be a scheme locally of finite type and $Y\subset X$ be an open subscheme or a closed subscheme of $X$. Then
$$\Sec_X(\dist(Y))\simeq\Gamma(Y)$$
where $\Gamma$ is the global section of schemes.
\end{cor}
\begin{proof}
Note that any open subschemes and closed subschemes of $X$ are locally of finite type so that there is a natural injective homomorphism $\Gamma(Y)\hookrightarrow\dist(Y)^*$. We may view $\Gamma(Y)$ as a subalgebra of $\dist(Y)^*$. Let $f:\spc(A)\rightarrow X$ be a morphism of schemes and $g:=\dist(f):A^\circ\rightarrow\dist(X)$ be a corresponding coalgebra homomorphism. Then $g^\dd(\dist(Y))=\dist(f^{-1}(Y))$ since $\dist$ preserves limits. The dual $g_{\dist(Y)}^*:\dist(Y)^*\rightarrow\dist(f^{-1}(Y))^*$ of $g_{\dist(Y)}$ introduced before \ref{secdef} restricts to an algebra homomorphism $f_Y:\Gamma(Y)\rightarrow\Gamma(f^{-1}(Y))$ induced by the morphism $f$ of schemes, i.e., the following commutes:
\[\begin{tikzcd}
	{\dist(Y)^*} & {\dist(f^{-1}(Y))^*} \\
	{\Gamma(Y)} & {\Gamma(f^{-1}(Y))}
	\arrow["{g^*_{\dist(Y)}}", from=1-1, to=1-2]
	\arrow[hook, from=2-1, to=1-1]
	\arrow["{f_Y}", from=2-1, to=2-2]
	\arrow[hook, from=2-2, to=1-2]
\end{tikzcd}\]
If $Y\subset X$ is open(resp. closed), then $A_{\dist(f^{-1}(Y))}=\Gamma(f^{-1}(Y))$ by \ref{dsec}(resp.by \ref{ac}). This shows that $g^*_{\dist(Y)}(\Gamma(Y))\subset \Gamma(f^{-1}(Y))$. Therefore $\Gamma(Y)$ is contained in $\Sec_X(\dist(Y))$. 

To show the opposite containment, take an element 
$\phi\in\Sec_X(\dist(Y))$ and an affine open covering $X=\bigcup_i U_i$. The dual map $j_i^*:\dist(Y)^*\rightarrow\dist(Y\cap U_i)^*$ of the inclusion $j_i:\dist(Y\cap U_i)\subset\dist(Y)$ restricts to an algebra homomorphism from $\Sec_X(\dist(Y))$ to $\Gamma(U_i)_{\dist(Y\cap U_i)}$. If $Y$ is open(resp. closed), then $\Gamma(U_i)_{\dist(Y\cap U_i)}=\Gamma(Y\cap U_i)$ by \ref{dsec}(resp.by \ref{ac}). Then $\{i_j^*(\phi)\in\Gamma(Y\cap U_i)\}$ can be glued to obtain a unique element $\phi'\in\Gamma(Y)$. If we view $\phi'$ as an element of $\Sec_X(\dist(Y))$, then $j^*_i(\phi')=j^*_i(\phi)$ for all $i$ by construction.
Since the algebra homomorphism
$$\dist(Y)^*\rightarrow \prod_i\dist(Y\cap U_i)^*$$
defined by $\psi\mapsto (j^*_i(\psi))_i$ is injective, we have $\phi=\phi'\in\Gamma(Y)$.
\end{proof}

%Let $i:Y\rightarrow X$ be a locally closed immersion, i.e. $i$ is a composition $i=j_o\circ j_c$ of a closed immersion $j_c$ and an open immersion $j_o$. Then the induced coalgebra homomorphism $\dist(i):\dist(Y)\rightarrow\dist(X)$ is injective. In this way $\dist(Y)$ can be seen as a subcoalgebra of $\dist(X)$.

A $k$-scheme $Y$ is called a locally closed subscheme of $X$ if it is a closed subscheme of an open subscheme of $X$. The next corollary immediately follows from the previous corollary.

\begin{cor}
    Let $X$ be a scheme locally of finite type and $Y\subset X$ be a locally closed subscheme. Then 
    $$\Sec_X(\dist(Y))\simeq\Gamma(Y).$$
\end{cor}

We have seen that the functor $(-)^\circ$ fully-faithfully embeds the category $\frfd^{op}$ into the category $\csu$ which gives natural bijections between hom-sets as in . The following shows that the functor $\emb$ embedds $\sch^{lf}$ into $\csu$ as well.

\begin{cor}
For a scheme $X$ locally of finite type, the ringed coalgebra $\emb(X)$ satisfies condition $(A)$. In particular, the functor $\emb$ corestricts to a functor $\sch^{lf}\rightarrow\csu$.
\end{cor}
\begin{proof}
Corollary \ref{ssec} implies that $\Gamma(\dist(X))=\O_X(X)$ and hence the coalgebra homomorphism $\dist(X)\rightarrow\Gamma(\dist(X))^\circ$ is the morphism obtained by applying $\emb$ to the natural morphism $X\rightarrow\spc(\O_X(X))$ of schemes. Hence $\dist(X)$ is an object of $\csu$.
\end{proof}

We close this section by showing that the functor $\emb$ is fully faithful when $k$ is algebraically closed.
%it is restricted to the full subcategory of $\rsp$ whose objects are integral schemes locally of finite type when $k$ is algebraically closed. 

We briefly discuss sober topological spaces and sobrification of topological spaces. Recall that a topological space is said to be sober if every irreducible closed subset has a unique generic point. The underlying spaces of schemes are known to be sober. For a topological space $X$, we denote by $S(X)$ the set of nonempty irreducible closed subsets of $X$. If $T\subset X$ is closed, then $S(T)\subset S(X)$. $S(X)$ can be endowed with a topology where the closed subsets are of the form $S(T)$ for some closed subset $T\subset X$. $S(X)$ is sober and this construction together with topological closure of images of continuous functions defines a functor from $\tops$ to the full subcategory $\sob$ of $\tops$ whose objects are sober spaces. This topological space $S(X)$ is known as the sobrification of $X$ and the functor is the left adjoint to the inclusion functor from $\sob$ to $\tops$.

\begin{lem}\label{top}
Let $X$ be a sober topological space such that the subset of closed points is dense. Then $X$ is naturally isomorphic to $S(\lv X\rv)$. 
\end{lem}
\begin{proof}
Since $X$ is sober, the map $x\mapsto \overline{\{x\}}$ gives an isomorphism from $X$ to $S(X)$. It is enough to show that the map $T\mapsto T\cap\lv X\rv$ is an isomorphism from $S(X)$ to $S(\lv X\rv)$. The map is bijective since the subset of closed points of $X$ is dense. Note that a closed subset $T'\subset X$ is irreducible if and only if $T'\cap\lv X\rv$ is. Therefore 
$S(T)\subset S(X)$ for some closed $T\subset X$ correspond to $\{T'\cap\lv X\rv| T'\in S(T)\}=S(T\cap\lv X\rv)\subset S(\lv X\rv)$.
\end{proof}

Next, we prove that a corestriction of a morphism of $\cs$ is again a morphism under certain assumptions. %For a ringed coalgebra $C$ and a subcoalgebra $D\subset C$, we denote by 

\begin{lem}\label{res}
Let $C_1, C_2, D_1$ and $D_2$ be ringed coalgebras satisfying (A), $f, g, h$ and $i$ are coalgebra homomorphisms making the following diagram commute: 
\[\begin{tikzcd}
	{C_1} & {C_2} \\
	{D_1} & {D_2}
	\arrow["g", from=1-1, to=1-2]
	\arrow["f", from=2-1, to=1-1]
	\arrow["h"', from=2-1, to=2-2]
	\arrow["i"', from=2-2, to=1-2]
\end{tikzcd}\]
If $f$ and $g$ are morphisms of ringed coalgebras and the composition 
$$j:D_2\overset{i}{\rightarrow} i(D_2)\hookrightarrow i(D_2)^{*\circ}\rightarrow\Sec_{C_2}(i(D_2))^\circ$$ 
of $i$ and the natural coalgebra homomorphisms is an isomorphism of ringed coalgebras(hence $i$ is injective), then $h$ is a morphism of ringed coalgebras.
\end{lem}
\begin{proof}
By commutativity of the diagram, the image of $D_1$ under $f$ is contained in $g^{-1}(i(D_2))$. Since $f(D_1)\subset C_1$ is a subcoaglebra, we obtain the following commutative diagram:
\[\begin{tikzcd}
	{g^\dd(i(D_2))} & {i(D_2)} \\
	{D_1} & {D_2}
	\arrow["g", from=1-1, to=1-2]
	\arrow["f", from=2-1, to=1-1]
	\arrow["h"', from=2-1, to=2-2]
	\arrow["i"', from=2-2, to=1-2]
\end{tikzcd}\]
Here we abuse notations to denote the coalgebra homomorphisms induced by $f$, $g$ and $i$.

We first check that the following diagram commutes:
\[\begin{tikzcd}
	{D_1} & {g^\dd(i(D_2))} \\
	{D_1^{*\circ}} & {g^\dd(i(D_2))^{*\circ}} \\
	{\Sec_{D_1}(D_1)^\circ} & {\Sec_{C_1}(g^\dd(i(D_2)))^\circ}
	\arrow["f", from=1-1, to=1-2]
	\arrow["nat.", hook, from=1-1, to=2-1]
	\arrow["nat.", hook, from=1-2, to=2-2]
	\arrow["f^{*\circ}", from=2-1, to=2-2]
	\arrow[from=2-1, to=3-1]
	\arrow[from=2-2, to=3-2]
	\arrow["{mor.}", from=3-1, to=3-2]
\end{tikzcd}\]
Here $nat.$ stands for the natural transformation $id_{\cog}\Rightarrow (-)^{*\circ}$. The bottom square is the diagram obtained by applying the functor $(-)^\circ:\cog\rightarrow\alg$ to the square
\[\begin{tikzcd}
	{g^\dd(i(D_2))^{*}} & {D_1^{*}} \\
	{\Sec_{C_1}(g^\dd(i(D_2)))} & {\Sec_{D_1}(D_1)}
	\arrow["{f^*}", from=1-1, to=1-2]
	\arrow[hook, from=2-1, to=1-1]
	\arrow[from=2-1, to=2-2]
	\arrow[hook, from=2-2, to=1-2]
\end{tikzcd}\]
which is commutative since $f$ is a morphism of ringed coalgebras.

Next, we check that the following diagram commutes:
\[\begin{tikzcd}
	{g^\dd(i(D_2))} & {i(D_2)} \\
	{g^\dd(i(D_2))^{*\circ}} & {i(D_2)^{*\circ}} \\
	{\Sec_{C_1}(g^\dd(i(D_2))^\circ} & {\Sec_{C_2}(i(D_2))^\circ}
	\arrow["{g}", from=1-1, to=1-2]
	\arrow["nat.", hook, from=1-1, to=2-1]
	\arrow["nat.", hook, from=1-2, to=2-2]
	\arrow["{g}^{*\circ}", from=2-1, to=2-2]
	\arrow[from=2-1, to=3-1]
	\arrow[from=2-2, to=3-2]
	\arrow["{mor.}", from=3-1, to=3-2]
\end{tikzcd}\]
Again, $nat.$ stands for the natural transformation $id_{\cog}\Rightarrow (-)^{*\circ}$. The bottom square is the diagram obtained by applying the functor $(-)^\circ:\cog\rightarrow\alg$ to the following square
\[\begin{tikzcd}
	{i(D_2)^{*}} & {g^\dd(i(D_2))^{*}} \\
	{\Sec_{C_2}(i(D_2))} & {\Sec_{C_1}(g^\dd(i(D_2))}
	\arrow["{g^*}", from=1-1, to=1-2]
	\arrow[hook, from=2-1, to=1-1]
	\arrow[from=2-1, to=2-2]
	\arrow[hook, from=2-2, to=1-2]
\end{tikzcd}\]
which is commutative since $g$ is a morphism of ringed coalgebras.

Combining all the diagrams, we have the following commutative diagram:
\[\begin{tikzcd}
	&& {D_2} \\
	{D_1} & {g^\dag(i(D_2))} & {i(D_2)} \\
	{\Sec_{D_1}(D_1)^\circ} & {\Sec_{C_1}(g^\dag(i(D_2)))^\circ} & {\Sec_{C_2}(i(D))^\circ}
	\arrow["i", hook, from=1-3, to=2-3]
	\arrow["j", curve={height=-20pt}, from=1-3, to=3-3]
	\arrow["h", from=2-1, to=1-3]
	\arrow["f"', from=2-1, to=2-2]
	\arrow["mor.", from=2-1, to=3-1]
	\arrow["{g}"', from=2-2, to=2-3]
	\arrow[from=2-2, to=3-2]
	\arrow[from=2-3, to=3-3]
	\arrow["mor.", from=3-1, to=3-2]
	\arrow["mor.", from=3-2, to=3-3]
\end{tikzcd}\]
Here the leftmost vertical arrow is a morphism of ringed coalgebras since $D_1$ satisfies condition $(A)$. Also, we have seen that the bottom arrows are morphisms of ringed coalgebras. Since $j$ is an isomorphism, $h$ is a morphism.
\end{proof}

The following fact is used to recover a morphism of schemes from a morphism of ringed in the proof of \ref{intfull}:

\begin{lem}\label{afmor}
    Let $A$ and $B$ be finitely generated algebras. Then any morphism $\spc(B)\rightarrow \spc(A)$ of ringed spaces is a morphism of schemes.
\end{lem}
\begin{proof}
    Let $X:=\spc(B), Y:=\spc(A)$ and $f:X\rightarrow Y$ be a morphism of ringed spaces. Let $\varphi:A\rightarrow B$ be the algebra homomorphism between the global sections induced by $f$. We show that the morphism $g=\spc(\varphi):X\rightarrow Y$ of schemes induced by $\varphi$ is equal to $f$ as morphisms of schemes. 
    
    The morphism $f$ induces  a morphism $f_\p:A_{f(\p)}\rightarrow B_\p$ between the stalks for each $\p$ and it makes the following commute:
\[\begin{tikzcd}
	A & B \\
	{A_{f(\p)}} & {B_\p}
	\arrow["\varphi", from=1-1, to=1-2]
	\arrow["l_{f(\p)}", from=1-1, to=2-1]
	\arrow["l_{\p}", from=1-2, to=2-2]
	\arrow["{f_\p}", from=2-1, to=2-2]
\end{tikzcd}\]
where $l_\p$ and $l_{f(\p)}$ stand for the canonical algebra homomorphisms of localizations. If $\p$ is a maximal ideal $\m\subset B$, then $B/\m$ is isomorphic to $k$ by the Nullstellensatz and so is $B_\m/\m B_\m$. Thus $f_\m^{-1}(\m B_\m)=f(\m)A_{f(\m)}$. 
By computing the kernel of the algebra homomorphism from $A$ to $k$, we have
$$\varphi^{-1}(\m)=\ker(p_\m\circ l_\m\circ\varphi)=\ker(p_\m\circ f_\m\circ l_{f(\m)})=A\cap f(\m)A_{f(\m)}=f(\m)$$
where $p_\m$ is the unique algebra homomorphism $B_\m\rightarrow k$.

Therefore $f$ and $g$ agree on the subset $\lv X\rv$ of closed points of $X$. The density of the subset $\lv X\rv\subset X$ and separation of affine schemes imply that
$f=g$ as continuous functions on $X$.

Again, for every $\p\in X$ we have the following commutative diagram
\[\begin{tikzcd}
	A & B \\
	{A_{\varphi^{-1}(\p)}} & {B_\p}
	\arrow["\varphi", from=1-1, to=1-2]
	\arrow["l_{f(\p)}", from=1-1, to=2-1]
	\arrow["l_{\p}", from=1-2, to=2-2]
	\arrow["{f_\p}", from=2-1, to=2-2]
\end{tikzcd}\]
It is a routine to check that $s\notin\p$ if $s\notin\varphi^{-1}(\p)$ for every $\p\in X$ and $s\in A$ and that
$$f_\p\left(\frac{a}{s}\right)=\frac{\varphi(a)}{\varphi(s)}.$$
Therefore, $f=g$ as morphisms of schemes.
\end{proof}

\begin{thm}\label{intfull}
%The restriction of $\emb$ to the full subcategory of integral schemes locally of finite type is
Let $k$ be an algebraically closed field. Then $\emb$ is fully-faithful, i.e., the map
$$\emb:\sch^{lf}(X, Y)\rightarrow\csu(\emb(X), \emb(Y))$$
is bijective for any schemes $X, Y$ locally of finite type over $k$.
\end{thm}
\begin{proof}
By corollary \ref{tfaith1}(or \ref{tfaith2}), it suffices to show the surjectivity. Let $g:\dist(X)\rightarrow \dist(Y)$ be a morphism of ringed coalgebras. Then it induces a continuous function $g|_{\lv X\rv}:\lv X\rv\rightarrow\lv Y\rv$ by \ref{com}. Since $X$ and $Y$ are sober and the subsets $\lv X\rv\subset X$ and $\lv Y\rv\subset Y$ are dense, it extends to a continuous function $f:X\rightarrow Y$ by \ref{top}. 

We show that $f$ extends to a morphism of ringed spaces. Let $U\subset Y$ be an affine open subscheme. The image of a pointed irreducible coalgebra under a coalgebra homomorphism is a pointed irreducible coalgebra(\cite[8.0.9]{S}), so for each $x\in \lv X\rv$ the image of $\O_{X, x}^\circ$ under $g$ lies in $\O_{Y, g(x)}^\circ$. Thus $g$ restricts to a coalgebra homomorphism $g|_{\dist(f^{-1}(U))}$ from $\dist(f^{-1}(U)))$ to $\dist(U)$ and we have the following commutative diagram:
\[\begin{tikzcd}
	{\dist(X)} & {\dist(Y)} \\
	{\dist(f^{-1}(U))} & {\dist(U)}
	\arrow["g", from=1-1, to=1-2]
	\arrow["incl.", hook, from=2-1, to=1-1]
	\arrow["{g|_{\dist(f^{-1}(U))}}"', from=2-1, to=2-2]
	\arrow["incl.", hook, from=2-2, to=1-2]
\end{tikzcd}\]
Since $U$ is affine, the natural coalgebra homomorphism $\dist(U)\rightarrow\Sec_Y(\dist(U))^\circ$ is an isomorphism of ringed coalgebras by \ref{ssec}. 

Therefore the restriction 
$$g|_{\dist(f^{-1}(U))}:\dist(f^{-1}(U))\rightarrow\dist(U)$$ 
is a morphism of ringed coalgebras by \ref{res}. Thus we obtain an algebra homomorphism 
$$f_U:=\Gamma(g|_{\dist(f^{-1}(U))}):\O_Y(U)\rightarrow\O_X(f^{-1}(U))$$ 
by \ref{dsec} as a restriction of the dual $g|_{\dist(f^{-1}(U))}^*$. For an open set $U\subset Y$ that is not necessarily affine, the glueing property of schemes induces an algebra homomorphism
$f_U:\O_Y(U)\rightarrow\O_X(f^{-1}(U))$. The functoriality of $\dist$ implies that the collection $\{f_U\}_U$ of homomorphisms is compatible with restrictions of sections. Therefore $f:X\rightarrow Y$ defines a morphism of ringed spaces. By \ref{afmor}, $f$ is indeed a morphism of schemes. 

Next, we show $\dist(f)=g$. By construction, $\dist(f)$ and $g$ agree on $\lv X\rv$. It suffices to show that they agree on $\O_{X, x}^\circ$ for all $x\in \lv X\rv$. Take an affine open neighborhood $V\subset f^{-1}(U)$ of $x$ and denote by $\dist(f)_V$ and by $g_V$ the coalgebra homomorphism $\dist(V)\rightarrow\dist(U)$ induced by $\dist(f)$ and $g$. Then the following diagram commutes:
\[\begin{tikzcd}
	& {\dist(f^{-1}(U))} \\
	& {\O_X(f^{-1}(U))^\circ} \\
	{\dist(V)=\O_X(V)^\circ} && {\dist(U)=\O_Y(U)^\circ}
	\arrow[from=1-2, to=2-2]
	\arrow["{g|_{\dist(f^{-1}(U))}}", from=1-2, to=3-3]
	\arrow["{f_U^\circ}"', dashed, from=2-2, to=3-3]
	\arrow[hook, from=3-1, to=1-2]
	\arrow["{\rho^{f^{-1}(U)\circ}_V}"', from=3-1, to=2-2]
	\arrow["{g_V}"', from=3-1, to=3-3]
\end{tikzcd}\]
In particular, the bottom triangle is commutative. Therefore $f^\circ_U$ and $g_V$ induce the same coalgebra homomorphism $\O_{X, x}^\circ\rightarrow\dist(U)$.
%$$\dist(f)_V=(\rho^{f^{-1}(U)}_V\circ f_U)^\circ=f^\circ_U\circ\rho^{f^{-1}(U)\circ}_V=g_V.$$ 
This shows that $g=\dist(f)$.
\end{proof}

Recall that $\aff$ denotes the full subcategory of $\alg$ whose objects are commutative finitely generated algebras. By combining \ref{com} and \ref{intfull}, we obtain the following:

\begin{cor}\label{ffcom}
Let $k$ be an algebraically closed field. The following diagram commutes up to natural isomorphism:
\[\begin{tikzcd}
	{\aff^{op}} && {\frfd^{op}} \\
	\sch^{lf} && \csu \\
	& \tops
	\arrow["incl.", hook, from=1-1, to=1-3]
	\arrow["\spc"', hook, from=1-1, to=2-1]
	\arrow["{(-)^\circ}", hook, from=1-3, to=2-3]
	\arrow[hook, "\emb", from=2-1, to=2-3]
	\arrow["{\lv-\rv}"', from=2-1, to=3-2]
	\arrow["pts", from=2-3, to=3-2]
\end{tikzcd}\]
Here $\lv-\rv$ is the functor that gives the subspace of closed points of the underlying space of schemes.
\end{cor}

\section{Modules over ringed coalgebras}

In the scheme theory, sheaves of modules over a scheme play a role of modules over a commutative ring. Sheaves on a ringed space $X$ are  contravariant functors from the category of the open sets of $X$ to the category $\cring$. When $X=\spc(R)$ for some commutative ring $R$, then every left $R$-module $M$ gives a sheaf $\tilde{M}$ on $X$ and the assignment $M\mapsto \tilde{M}$ defines a fully-faithful functor.

In this section, we introduce modules over a ringed coalgebra. If $C=(C, Q, \Sec)$ is a ringed coalgebra, a module over $C$ is a functor from $P_C^{op}$ to the category $\vect$ satisfying some compatibility with the section $\Sec:P^{op}_C\rightarrow\rfd$ on $C$. The main ingredient of the construction is the comodules over the subcoalgebras of the underlying coalgebra $C$. We discuss modules over ringed coalgebras in a way similar to the constructions used for ringed coalgebras; subsections \ref{mod0}, \ref{mod1}, \ref{mod2} and \ref{mod3} correspond to subsections \ref{cog0}, \ref{cog1}, \ref{cog2} and \ref{cog3} respectively. We first briefly review its definition and some category theoretic properties in the first part \ref{mod0} of this section. By the same technique we used in the previous section, we will construct a module $\ha{M}$ over a ringed coalgebra $A^\circ$ using a left $A$-module $M$ in the second subsection \ref{mod1}. The assignment $M\mapsto \ha{M}$ defines a faithful functor and it is also full on the full subcategory of finitely generated $A$-modules. In subsection \ref{mod2}, we further apply the construction to $\O_X$-modules on a  scheme $X$ locally of finite type. Finally, in the last part \ref{mod3} of this section, we fully-faithfully embed the category of coherent sheaves on $X$ into the category of modules over the ringed coalgebra $\emb(X)$ when $X$ is separated. Except of \ref{mod3}, the base field $k$ is not necessarily algebraically closed.

\subsection{The finite dual comodules}\label{mod0}

Coalgebras are the underlying structure of ringed coalgebras. On the other hand, a module $\rcm$ over a ringed coalgebra $C$ is a functor $\rcm:P^{op}_C\rightarrow\vect$ that respects the action of $\Sec_C:P^{op}\rightarrow\rfd$, however, we will use comodules to define a module over a ringed coalgebra. We review a definition of comodules in this subsection.

Let $C$ be a coalgebra. A left comodule $M=(M, \rho)$ over $C$ is a pair of a vector space $M$ and a linear map $\rho:M\rightarrow C\otimes M$ that makes the diagram on the left hand side commute. A comodule homomorphism $f:(M_1, \rho_1)\rightarrow(M_2, \rho_2)$ is a linear map $f:M_1\rightarrow M_2$ that makes the diagram on the RHS commute.
\[\begin{tikzcd}
	M && {C\otimes M} && {M_1} && {M_2} \\
	\\
	{C\otimes M} && {C\otimes C\otimes M} && {C\otimes M_1} && {C\otimes M_2}
	\arrow["\rho", from=1-1, to=1-3]
	\arrow["\rho"', from=1-1, to=3-1]
	\arrow["{id_C\otimes\rho}", from=1-3, to=3-3]
	\arrow["f", from=1-5, to=1-7]
	\arrow["{\rho_1}"', from=1-5, to=3-5]
	\arrow["{\rho_2}", from=1-7, to=3-7]
	\arrow["{\Delta\otimes id_M}"', from=3-1, to=3-3]
	\arrow["{id_C\otimes f}"', from=3-5, to=3-7]
\end{tikzcd}\]
Through out this paper, every comodule over a coalgebra will be a left comodule unless stated otherwise. The category of left comodules over $C$ and comodule homomorphisms will be denoted by ${}_C\cod$.

For an algebra $A$, denote by ${}_A\modu$(resp.$\modu_A$) the category of left(resp. right) $A$-modules and their homomorphisms. If $A$ is an RFD algebra, then there is an adjunction between $_{A^\circ}\cod$ and ${}_{A}\modu^{op}$ similar to the one between and $\cog$ and $\alg^{op}$ which we describe below. For details, see \cite[Section 1.7]{T2}. 

Let $M=(M, \rho)$ be a left comodule over $A^\circ$. Then the dual space $M^*$ has a left $A^{\circ*}$-module structure defined by
$$f\cdot g(x):=(f\otimes g)(\rho(x))$$
for all $f\in A^{\circ*}$, $g\in M^*$ and $x\in M$. Since $A$ is a subalgebra of $A^{\circ*}$, $M^*$ is an $A$-module. This construction defines a functor $(-)^*:{}_{A^\circ}\cod\rightarrow {}_{A}\modu^{op}$.

Next, we describe a left adjoint of $(-)^*$. The finite dual $M^\circ$ of a left $A$-module $M$ is defined by setting
$$M^\circ=\{\phi\in M^*|\exists N\subset\ker\phi\ s.t.\ \text{$N$ is a left $A$-submodule of finite codimension}\}.$$ 
Below we show that $M^\circ$ has a left $A^\circ$-comodule structure. 

Recall that a right $A$-module $M$ is said to be locally finite if the subspace $xA\subset M$ is of finite dimension for every $x\in M$. Let $\modu_A^{lf}$ be the full subcategory of $\modu_A$ whose objects are locally finite modules. This is indeed a reflective subcategory of $\modu_A$: the left adjoint of the inclusion functor is given by
$$M^{lf}:=\{x\in M|xA\ \text{is of finite dimension}\}.$$
By symmetric argument in the lemma in \cite[1.7.1]{T2}, there is an equivalence of categories between $\modu_A^{lf}$ and ${}_{A^\circ}\cod$. This equivalence of categories assigns to every $A^\circ$-comodule $N$ the underlying vector space $N$ with the right $A$-action defined by
$$x\cdot a:=\sum_i\phi_i(a)x_i$$
for all $x\in N$ and $a\in A$ where $\rho(x)=\sum_i \phi_i\otimes x_i$ where $\phi_i\in A^\circ$ and $x_i\in N$. This action turns $N$ into a locally finite right $A$-module. All in all, we have the following adjunctions:
\[\begin{tikzcd}
	{(*)\ \ \ {}_A\modu^{op}} & {\modu_A} & {\modu^{lf}_A} & {{}_{A^\circ}\cod}
	\arrow[""{name=0, anchor=center, inner sep=0}, "{(-)^T}"', shift right=2, from=1-1, to=1-2]
	\arrow[""{name=1, anchor=center, inner sep=0}, "{(-)^T}"', shift right=2, from=1-2, to=1-1]
	\arrow[""{name=2, anchor=center, inner sep=0}, "{(-)^{lf}}"', shift right=2, from=1-2, to=1-3]
	\arrow[""{name=3, anchor=center, inner sep=0}, shift right=2, hook, from=1-3, to=1-2]
	\arrow["\simeq"{description}, draw=none, from=1-3, to=1-4]
	\arrow["\dashv"{anchor=center, rotate=90}, draw=none, from=0, to=1]
	\arrow["\dashv"{anchor=center, rotate=90}, draw=none, from=2, to=3]
\end{tikzcd}\]
 The functors on the left are the ones induced by transposed $A$-module structures: if $M$ is a left(resp.right) $A$-module, then $M^*$ can be endowed with a right(resp.left) $A$-module structure by
$$(\phi\cdot a)(x):=\phi(ax)(\text{resp.}(a\cdot \phi)(x):=\phi(xa))$$
where $\phi\in M^*$ and $a\in A$.  There is a natural commutative diagram
\[\begin{tikzcd}
	{\vect(M, N^*)} & {\vect(N, M^*)} \\
	{{}_A\modu(M, N^*)} & {\modu_{A}(N, M^*)}
	\arrow["\sim", from=1-1, to=1-2]
	\arrow[hook, from=2-1, to=1-1]
	\arrow["\sim", from=2-1, to=2-2]
	\arrow[hook, from=2-2, to=1-2]
\end{tikzcd}\]
%$${}_A\modu(M, N^*)\simeq\modu_A(N, M^*)$$
where $M$ is a left $A$-module and $N$ is a right $A$-module. 

In the above diagram $(*)$, the functor ${}_{A^\circ}\cod\rightarrow{}_A\modu^{op}$ is nothing but the dual module functor $(-)^*$. Thus the composition $(-)^{*lf}:{}_A\modu^{op}\rightarrow{}_{A^\circ}\cod$ is the left adjoint of $(-)^*$. The following statement follows from the same argument in the proof of \cite[1.7.2]{T2}.
\begin{lem}
Let $M$ be a left $A$-module. Then $M^{*lf}=M^\circ$ as subspaces of the dual space $M^*$.
\end{lem}

Therefore the subspace $M^\circ$ has a left $A^\circ$-comodule structure in such a way that if $\phi\in M^\circ$, there exist collections of finitely many elements $\{\psi_i\}\subset A^\circ, \{\phi_i\}\subset M^\circ$ such that
$$\phi(ax)=\sum_i\psi_i(a)\phi_i(x)$$
for all $a\in A$ and $x\in M$. We will call $M^\circ$ the \emph{finite dual comodule} of $M$. The construction of $M^\circ$ defines a functor $(-)^\circ:{}_A\modu^{op}\rightarrow{}_{A^\circ}\cod$. Furthermore, it is a left adjunction of $(-)^*:{}_{A^\circ}\cod\rightarrow {}_{A}\modu^{op}$. 

%This adjunction can be given as a composite of adjunctions between subcategories of $\modu_A$. 

%: one can assign to a locally finite module $M$ a comodule
%and every left $A^\circ$-comodule $N$ has a right $A$-module structure given by
%$$xa:=\sum_i \phi_i(a)x_i$$
%where the coaction $\rho$ is given by $\rho(x)=\sum_i \phi_i\otimes x_i, \phi_i\in A^\circ, x_i\in N$. This action makes $N$ a locally finite module.

%Together with the adjunction between $\modu_A$ and ${}_{A^\circ}\cod$ in \cite{T2}, we have the following diagram:

A left $A$-module $M$ will be called \emph{residually finite-dimensional}(RFD) if the intersection of all left submodules of finite codimension is $0$. This is equivalent to $M$ being a subdirect product of finite dimensional left $A$-modules. Note that $A$ is RFD as a left $A$-module if $A$ is RFD as an algebra. Indeed, if $I\subset A$ is a left ideal of finite codimension, then the annihilator
$$\ann(A/J)=\{x\in A|xA\subset J\}\subset J$$
of a finite dimensional by the argument in \ref{frfd}. 

We prepare the following lemma for later use.
\begin{lem}\label{modd}
Let $A$ be a commutative algebra. Then there is an isomorphism
$$M^\circ\simeq\bigoplus_{\m\in\lv\spc(A)\rv} M^\circ_\m$$
of left comodules over $A^\circ$ for every finitely generated $A$-module $M$. 
\end{lem}
\begin{proof}
Since ${}_{A^\circ}\cod$ and $\modu_{A}^{lf}$ are equivalent, it suffices to show that $M^\circ$ and $\bigoplus_{\m\in\lv\spc(A)\rv} M^\circ_\m$ are isomorphic as right $A$-modules. The algebra $A$ is commutative, so it suffices to show the isomorphy as left $A$-modules.

The \cite[1.7.5]{T2} states that there is an isomorphism
$$M^\circ\simeq{}_A\modu(M, A^\circ)$$
of left $A$-modules. Then
$$M^\circ\simeq {}_A\modu(M, A^\circ)\simeq{}_A\modu(M, \bigoplus_{\m\in\lv\spc(A)\rv} A^\circ_\m)$$
by \cite[2.1.7]{T2}. Since $M$ is finitely generated, the direct sum on the right hand side commutes with the hom-functor ${}_A\modu(M, -):{}_A\modu\rightarrow{}_A\modu$:
$${}_A\modu(M, \bigoplus_{\m\in\lv\spc(A)\rv}A^\circ_\m)\simeq\bigoplus_{\m\in\lv\spc(A)\rv}{}_A\modu(M, A^\circ_\m).$$ 
Therefore
$$\bigoplus_{\m\in\lv\spc(A)\rv}{}_A\modu(M, A^\circ_\m)\simeq\bigoplus_{\m\in\lv\spc(A)\rv}{}_{A_\m}\modu(M_\m, A^\circ_\m)\simeq\bigoplus_{\m\in\lv\spc(A)\rv} M^\circ_\m$$
as left $A$-modules.
\end{proof}

\subsection{Modules over  $A^\circ$}\label{mod1}

Every left module $M$ over a commutative ring $R$ produces a sheaf $\tilde{M}$ of modules over the spectrum $\spc(R)$ of the ring. 
The construction of $\tilde{M}$ is a generalization of that of the structure sheaf on $\spc(R)$. This gives a fully-faithful functor from the category ${}_R\modu$ to the category of quasi-coherent sheaves on $\spc(R)$ whose left adjoint is the global section functor. Similarly, every left module $M$ over an RFD algebra $A$ gives a module $\ha{M}$ over the ringed coalgebra $A^\circ$. We show that this gives a fully faithful functor to a certain full subcategory of ${}_{A^\circ}\modu$ and a theorem similar to \ref{algadj} holds.

\begin{de}
Let $C=(C, Q, \Sec_C)$ be a ringed coalgebra. We say a functor $\rcm:P_C^{op}\rightarrow\vect$ is a \emph{module over $C$} if for every subcoalgebra $D\subset C$, $\shm(D)$ is a RFD left $\Sec(D)$- module and the following diagram commutes for all subcoalgebras $D'\subset D\subset C$:
\[\begin{tikzcd}
	{\Sec(D)\otimes\rcm(D)} & {\Sec(D')\otimes\rcm(D')} \\
	{\rcm(D)} & {\rcm(D')}
	\arrow["{\rho^D_{D'}\otimes \lambda^D_{D'}}", from=1-1, to=1-2]
	\arrow["\mu_D", from=1-1, to=2-1]
	\arrow["\mu_{D'}", from=1-2, to=2-2]
	\arrow["\lambda^D_{D'}"', from=2-1, to=2-2]
\end{tikzcd}\]
where $\rho^D_{D'}$ and $\lambda^D_{D'}$ are restrictions defined in definition \ref{rc} and $\mu_D, \mu_{D'}$ are the actions of the algebras over the modules. 
If $\rcm_1, \rcm_2$ are modules over the ringed coalgebra $C$, then a morphism $\rcm_1\rightarrow\rcm_2$ is a natural transformation $\{f_D\}$ that it compatible with the action of $\Sec$, i.e., the following diagram commutes for all subcoalgebras $D\subset C$:
\[\begin{tikzcd}
	{\Sec(D)\otimes\rcm_1(D)} & {\Sec(D)\otimes\rcm_2(D)} \\
	{\rcm_1(D)} & {\rcm_2(D)}
	\arrow["{\id\otimes f_D}", from=1-1, to=1-2]
	\arrow["\mu^1_D", from=1-1, to=2-1]
	\arrow["\mu^2_D", from=1-2, to=2-2]
	\arrow["f_D"', from=2-1, to=2-2]
\end{tikzcd}\]
where $\mu^1_D$ and $\mu^2_D$ are actions of the algebra $\Sec(D)$ on the modules $\rcm_1(D)$ and $\rcm_2(D)$.

The category of modules over $C$ and their morphisms will be denoted by ${}_C\modu$.
\end{de}

We abuse notation to define the global section functor on modules over ringed coalgebras.

\begin{de}
Let $C=(C, Q, \Sec)$ be a ringed coalgebra.
    We define a functor $\Gamma:{}_{C}\modu\rightarrow _{\Sec(C)}\modu$ by
    $$\Gamma(\rcm):=\rcm(C)$$
    for every module $\rcm$ over $C$. This functor will be called the global section functor.
\end{de}

The rest of this subsection is devoted to construct a module $\ha{M}$ over a ringed coalgebra $A^\circ$ associated to a left $A$-module $C$. The construction uses cotensor product of comodules. Let $C$ be a coalgebra, $D\subset C$ be a subcoalgebra, and $N$ be a left comodule over $C$. The \emph{cotensor product} of $D$ and $N$ is given by
$$D\sq N:=\{x\in N|\rho(x)\in D\otimes N\}$$
where $\rho$ is the coaction of $C$ on $N$. This is a left $D$-comodule. $D\sq(-)$ commutes with direct sum: $D\sq(\bigoplus_i N_i)=\bigoplus_i(D\sq N_i)$. For details, see \cite{T3}.

Now let $A$ be a fully RFD algebra and $M$ be a finitely generated left $A$-module and $C\subset A^\circ$ be a subcoalgebra. Note that the natural map $M\rightarrow M^{\circ*}$ is injective since $A$ is a fully RFD algebra. The finite dual $M^\circ$ is a left $A^\circ$-comodule and hence $C\sq M^\circ$ is a left $C$-comodule. This induces the transposed right action of $C^*$ on $C\sq M^\circ$. Thus $(C\sq M^\circ)^*$ has left $C^*$-module structure. Recall from \ref{sec} that
$$\F_C:=\{I\subset A|I\ \text{is a left ideal},\ Z^\circ(I)\cap C=0\}.$$
We define
$$T^M_C:=\{x\in(C\sq M^\circ)^*|\exists I\in\F_C\ a_C(I)x\subset m^M_C(M)\}$$
where $a_C:A\rightarrow C^*$ is the algebra homomorphism defined in $\ref{sec}$ and $m^M_C$ is the composition of the natural module homomorphisms $M\hookrightarrow M^{\circ*}$ and $M^{\circ*}\twoheadrightarrow (C\sq M)^*$. We denote by $M_C$ the left $A_C$-submodule of $(C\sq M^\circ)^*$ generated by $T^M_C$:
$$M_C:=\Sec_X(C)T^M_C\subset (C\sq M^\circ)^*.$$
By an argument similar to \ref{ac}, we can show that if $C=Z^\circ(I)$ for some ideal $I\subset A$, then $\ha{M}(C)=M/IM$. Also, the argument in \ref{ac2} shows that if every ideal in $\F_C$ is two-sided(e.g. $A$ is commutative), then $T^M_C=M_C$. 

For subcoalgebras $A^\circ\supset C_1\supset C_2$, we have $C_2\sq M^\circ\subset C_1\sq M^\circ$. The dual map $(C_1\sq M^\circ)^*\rightarrow (C_2\sq M^\circ)^*$ of the inclusion induces a map $T^M_{C_1}\rightarrow T^M_{C_2}$ and hence a linear map $\lambda^{C_1}_{C_2}:M_{C_1}\rightarrow M_{C_2}$ that makes the following diagram commute:
\[\begin{tikzcd}
	{A_{C_1}\otimes M_{C_1}} & {A_{C_2}\otimes M_{C_2}} \\
	{M_{C_1}} & {M_{C_2}}
	\arrow["{\rho^{C_1}_{C_2}\otimes\lambda^{C_1}_{C_2}}"', from=1-1, to=1-2]
	\arrow[from=1-1, to=2-1]
	\arrow[from=1-2, to=2-2]
	\arrow["\lambda^{C_1}_{C_2}"', from=2-1, to=2-2]
\end{tikzcd}\]

\begin{de}
Let $A$ be a fully RFD algebra and $M$ be a finitely generated left $A$-module. We define a module $\ha{M}\in{}_{A^\circ}\modu$ associated to $M$ by
$$\ha{M}(C):=M_C$$
for every subcoalgebra $C\subset A^\circ$. 
\end{de}

Let $A$ be a fully RFD algebra. We show that the assignment $M\mapsto\ha{M}$ is functorial. For a homomorphism $f:M\rightarrow M'$ of left $A$-modules and a subcoalgebra $C\subset A^\circ$, we set
$$f^C:=\id_C\sq f^\circ:C\sq M'^{\circ}\rightarrow C\sq M^\circ.$$
Note that this is a morphism of left $C$-comodules.

\begin{lem}
    Let $A$ be a fully RFD algebra, $f:M\rightarrow M'$ be a homorphism of left $A$-modules and $C\subset A^\circ$ be a subcoalgebra. Then the dual $$f^{C*}:(C\sq M^\circ)^*\rightarrow (C\sq M'^\circ)^*$$
    restricts to a map from $T^M_C$ to $T^{M'}_C$ and hence a module homomorphism from $M_C$ to $M'_C$.
    \end{lem}
\begin{proof}
Let $x\in T^M_C$. Then $x\in(C\sq M^\circ)^*$ satisfies $a_C(I)x\subset m^M_C(M)$ for some $I\in\F_C$. Since $f^{C*}$ is a dual of a left $C$-comodule homomorphism, it is a left $C^*$-module homomorphism. Here $a_C$ is an algebra homomorphism from $A$ to $C^*$, so $a_C(I)\subset C^*$. Therefore 
$$a_C(I)f^{C*}(x)=f^{C*}(a_C(I)x)\subset f^{C*}(m^M_C(M)).$$
On the other hand, the following diagram commutes:
\[\begin{tikzcd}
	M & {M'} \\
	{(C\sq M^\circ)^*} & {(C\sq M'^\circ)^*}
	\arrow["f", from=1-1, to=1-2]
	\arrow["{m^M_C}"', from=1-1, to=2-1]
	\arrow["{m^{M'}_C}", from=1-2, to=2-2]
	\arrow["{f^{C*}}"', from=2-1, to=2-2]
\end{tikzcd}\]
Therefore we have
$$f^{C*}(m^M_C(M))\subset m^{M'}_C(f(M))\subset m^{M'}_C(M').$$ 
Thus 
$$a_C(I)f^{C*}(x)\subset m^{M'}_C(M')$$
so that $f^{C*}(x)\in T^{M'}_C$.
\end{proof}
 %Every $A$-module homomorphism $M\rightarrow M'$ induces a homomorphism $M_C\rightarrow M'_C$ of $A_C$-modules. Indeed, a $A$-module homomorphism $f:M_1\rightarrow M_2$ naturally induces a homomorphism
%$$f'=(C\sq f^\circ)^*:(C\sq M^\circ_1)^*\rightarrow (C\sq M^\circ_2)^*$$ 
%of $C^*$-modules for every subcoalgebra $C\subset A^\circ$. 

\begin{de}
    Let $A$ be a fully RFD algebra. The assignment that sends a finitely generated left $A$-module $M$ to the module $\ha{M}$ over the ringed coalgebra $A^\circ$ defines a functor $\ha{(-)}:{}_A\modu_{fg}\rightarrow{}_{A^\circ}\modu$.
\end{de}

\begin{cor}
Let $A$ be a fully RFD algebra. Then $\Gamma\circ (-)^\circ:{}_A\modu_{fg}\rightarrow{}_A\modu_{fg}$ is equal to the identity functor $\id_{{}_A\modu_{fg}}$ on ${}_A\modu_{fg}$.
\end{cor}
\begin{proof}
The equality $\Gamma(M^\circ)=M$ follows from \ref{ac}:
$$\Gamma(M^\circ)=M_{A^\circ}=m^M_{A^\circ}(A^\circ)=M.$$
The equality $\Gamma\circ (-)^\circ=id_{\frfd^{op}}$ on the morphisms of $\frfd^{op}$ follows from the fact that the diagram
\[\begin{tikzcd}
	{M^{\circ*}} & {M'^{\circ*}} \\
	M & M'
	\arrow["{f^{\circ*}}", from=1-1, to=1-2]
	\arrow[hook, from=2-1, to=1-1]
	\arrow["f", from=2-1, to=2-2]
	\arrow[hook, from=2-2, to=1-2]
\end{tikzcd}\]
commutes for every $A$-module homomorphism $f:M\rightarrow M'$.
\end{proof}

%Hence the construction  There is also a functor $\Gamma:{}_{A^\circ}\modu\rightarrow{}_A\modu_{fg}$ defined by
%$$\Gamma(\shm):=\shm(A^\circ).$$

%We denote this modules over $A^\circ$ by $M^\circ$. 

For every module $\rcm$ over a ringed coalgebra $C$, the inclusion $C'\subset C$ induces $\lambda^C_{C'}:\rcm(C)\rightarrow\rcm(C')$. The finite dual $(\lambda^C_{C'})^\circ:\rcm(C')^\circ\rightarrow\rcm(C)^\circ$ corestricts to $l^C_{C'}:\rcm(C')^\circ\rightarrow C'\sq\rcm(C)^\circ$. 

Next, we define a full subcategory ${}_C\modb$ of ${}_C\modu$ whose objects $\rcm$ satisfying the following property:
\begin{enumerate}
\item[(B)] The dual $(l^C_{C'})^{\circ*}:\rcm(C')^\circ\rightarrow C'\sq\rcm(C)^\circ$ restricts to an algebra homomorphism $\iota_C:\Gamma(\rcm)_{C'}\rightarrow\rcm(C')$ for all subcoalgebra $C'\subset C$.
\end{enumerate}
%$$ \text{(B)\ The dual $(l^C_{C'})^{\circ*}$ restricts to a map $\iota_C:\Gamma(\shm)_{C'}\rightarrow\shm(C')$ for all subcoalgebra $C'\subset C$.}$$
In other words, the following diagram commutes:
\[\begin{tikzcd}
	{(C'\sq\rcm(C)^\circ)^*} & {\rcm(C')^{\circ*}} \\
	{\rcm(C)_{C'}} & {\rcm(C')}
	\arrow["(l^C_{C'})^{\circ*}", from=1-1, to=1-2]
	\arrow["{incl.}"', hook, from=2-1, to=1-1]
	\arrow["{\iota_{C'}}"', dashed, from=2-1, to=2-2]
	\arrow["{incl.}"', hook, from=2-2, to=1-2]
\end{tikzcd}\]

Note that if $C=A^\circ$ for some fully RFD algebra, then the module of the form $\ha{M}$ for some left $A$-module $M$ satisfies (B). Indeed, $\rcm(A^\circ)_{C'}=M_{C'}$ and the dual map $(C'\sq M^\circ)^*\rightarrow M_{C'}^{\circ*}$ restricts to the identity on $M_{C'}=\rcm(C')$.

\begin{thm}
Let $A$ be a fully RFD algebra, $M$ be a finitely generated $A$-module and $\rcm$ be a $A^\circ$-module satisfying (B). There is a natural bijective correspondence
$${}_{A^\circ}\modb(\ha{M}, \rcm)\simeq {}_A\modu(M, \Gamma(\rcm))$$
given by  $f\mapsto \Gamma(f)$. 
\end{thm}
\begin{proof}
We show that the assignment $g\mapsto \iota\circ\ha{g}$ gives the inverse of $f\mapsto \Gamma(f)$. It is immediate to see that $g=\Gamma(\iota\circ\ha{g})$. The equality $f=\iota\circ\ha{\Gamma(f)}$ follows from the following commutative diagram:
\[\begin{tikzcd}
	{M^{\circ*}} &&& {\Gamma(\rcm)^{\circ*}} \\
	\\
	& {(C\sq M^\circ)^*} &&& {(C\sq\Gamma(\rcm)^\circ)^*} \\
	M &&& {\Gamma(\rcm)} \\
	&&& {M_{C}^{\circ*}} &&& {\rcm(C)^{\circ*}} \\
	& {M_{C}} &&& {\Gamma(\rcm)_{C}} \\
	&&& {M_{C}} &&& {\rcm(C)}
	\arrow["{\Gamma(f)^{\circ*}}", from=1-1, to=1-4]
	\arrow["{can.}"', two heads, from=1-1, to=3-2]
	\arrow["{\lambda_{M}^{\circ*}}"{pos=0.3}, curve={height=-24pt}, from=1-1, to=5-4]
	\arrow["{can.}"', two heads, from=1-4, to=3-5]
	\arrow["{\lambda_\rcm^{\circ*}}", curve={height=-24pt}, from=1-4, to=5-7]
	\arrow["{(id_C\sq \Gamma(f)^\circ)^*}"{pos=0.4}, from=3-2, to=3-5]
	\arrow["l_M^{\circ*}"'{pos=0.7}, from=3-2, to=5-4]
	\arrow["l_\shm^{\circ*}"'{pos=0.3}, from=3-5, to=5-7]
	\arrow["{incl.}", hook, from=4-1, to=1-1]
	\arrow["{\Gamma(f)}"'{pos=0.2}, from=4-1, to=4-4]
	\arrow["{\lambda_M}"', from=4-1, to=6-2]
	\arrow["{incl.}"{pos=0.7}, hook, from=4-4, to=1-4]
	\arrow["\lambda_{\Gamma(\rcm)}"'{pos=0.7}, from=4-4, to=6-5]
	\arrow["{\lambda_\rcm}", from=4-4, to=7-7]
	\arrow["{f_{C}^{\circ*}}"{pos=0.7}, from=5-4, to=5-7]
	\arrow["{incl.}"'{pos=0.3}, hook, from=6-2, to=3-2]
	\arrow["{\Gamma(f)_{C}}"{pos=0.3}, from=6-2, to=6-5]
	\arrow["id"', from=6-2, to=7-4]
	\arrow["{incl.}", hook, from=6-5, to=3-5]
	\arrow["{\iota_C}", from=6-5, to=7-7]
	\arrow["{incl.}"{pos=0.3}, hook, from=7-4, to=5-4]
	\arrow["{f_{C}}"'{pos=0.7}, from=7-4, to=7-7]
	\arrow["{incl.}", hook, from=7-7, to=5-7]
\end{tikzcd}\]
Here $\lambda_M$ and $\lambda_\rcm$ are restrictions of $\ha{M}$ and $\rcm$, $l_M^{\circ*}:(C\sq M^\circ)^*)\rightarrow M_C^{\circ*}$ and $l_\rcm^{\circ*}:(C\sq \Gamma(\rcm)^\circ)^*)\rightarrow \rcm(C)^{\circ*}$ are natural linear maps.
\end{proof}

\begin{cor}\label{embmod}
For any fully RFD algebra $A$, the functor $\ha{(-)}:{}_A\modu_{fg}\rightarrow{}_{A^\circ}\modb$ is fully-faithful, i.e.,
$$\ha{(-)}:{}_A\modu_{fg}(M_1, M_2)\rightarrow{}_{A^\circ}\modb(\ha{M}_1, \ha{M}_2)$$
is bijective.
\end{cor}
\begin{proof}
%Since $\Gamma\circ\ha{(-)}=id$, the map is injective.
It follows from the natural isomorphisms
$${}_{A^\circ}\modb(\ha{M}_1, \ha{M}_2)\simeq{}_A\modu(M_1, \Gamma(\ha{M}_2))\simeq{}_A\modu(M_1, M_2)$$
for any finitely generated $A$-modules $M_1$ and $M_2$.
\end{proof}

\subsection{Modules over $\dist(X)$}\label{mod2}
In this subsection, we discuss modules over the ringed coalgebra $\dist(X)$ associated to a scheme $X$ locally of finite type. The goal of this section is to associate a module over $\dist(X)$ to every coherent sheaf over $X$ by following the technique used in \ref{secdef}. We first introduce underlying comodules of a coherent sheaf, which is an analogue of underlying coalgebras of schemes.

\begin{de}
    Let $X$ be a scheme locally of finite type and $\shm$ be a coherent sheaf on $X$. We define a \emph{underlying comodule} $\dist(\shm)$ of $\shm$ by
$$\dist(\shm):=\bigoplus_{x\in\lv X\rv} \shm_{X, x}^\circ.$$
\end{de}

If $X=\spc(A)$ for some finitely generated algebra $A$ and $M$ is a finitely generated left $A$-module, then 
$$\dist(\til{M})=\bigoplus_{x\in\lv X\rv} \til{M}_x^\circ=\bigoplus_{x\in\lv \spc(A)\rv} M_x^\circ\simeq M^\circ$$ 
by \ref{modd} where $\til{M}$ is the sheaf associated to $M$. Note that a morphism $f:\shm_1\rightarrow\shm_2$ of sheaves induces a comodule morphism $\dist(f):\dist(\shm_2)\rightarrow\dist(\shm_1)$. Therefore the construction gives a functor $\dist:\coh_X^{op}\rightarrow{}_{\dist(X)}\modu$.
%$$M^\circ\simeq A^\circ\sq M^\circ=(\bigoplus A^\circ_\m)\sq M^\circ\simeq\bigoplus(A^\circ_\m\sq M^\circ)\simeq\bigoplus M^\circ_\m$$

%Let $\{U_i\}$ be the set of all affine open subschemes of $X$ whose global sections are RFD. 
Let $X$ be a scheme locally of finite type and $\shm$ be a coherent sheaf over $X$. Now we define a module $\ha{\shm}:P^{op}_{\dist(X)}
\rightarrow \vect$ 
over the ringed coalgebra $\dist(X)$ associated to a coherent sheaf $\shm$. Let $C\subset \dist(X)$ be a subcoalgebra and $U\subset X$ be an affine open subset. Then we have an inclusion $C\cap\dist(U)\subset C$
of subcoalgebras of $\dist(X)$. Then we obtain an inclusion of $C$-comodules:
$$i_C:(C\cap\dist(U))\sq\dist(\shm)\hookrightarrow C\sq\dist(\shm).$$
The comodule on the left hand side can be written as a subcomodule of the finite dual comodule associated to some 
finitely generated $\O_{X}(U)$-module.
Indeed, $C\cap \dist(U)$ is a subcoalgebra of $\dist(U)$ so that we have
$$(C\cap\dist(U))\sq\dist(\shm)=(C\cap\dist(U))\sq\dist(\shm|_U).$$
The restriction of $\shm$ to $U$ is the sheaf $\til{M}$ of modules over $U$ associated to some finitely generated $\O_X(U)$-module $M$. Therefore
$$(C\cap\dist(U))\sq\dist(\shm|_U)=
(C\cap\dist(U))\sq\dist(\til{M})$$
and we obtain an inclusion
$$i_C:(C\cap\dist(U))\sq\dist(\til{M})\subset C\sq\dist(\shm).$$
The dual map
$$i_C^*:(C\sq\dist(\shm))^*\twoheadrightarrow((C\cap\dist(U))\sq\dist(\til{M}))^*$$
is a homomorphism of left $C^*$-modules and can be thought of as a restriction of sections. Every element of $\shm(U)$ is locally an element of $\til{M}(U)$, i,e., restricted to an element of $\til{M}(U)$ on every affine open subset $U\subset X$. Similarly, we want every element in $\ha{\shm}(C)$ to locally be an element of $\ha{M}(C\cap\dist(U))$. In other words, there should be a map
$\ha{\shm}(C)\rightarrow \ha{M}(C\cap\dist(U))$ that make the following diagram commute:
\[\begin{tikzcd}
	{(C\sq\dist(\shm))^*} & {((C\cap\dist(U))\sq\dist(\til{M}))^*} \\
	{\ha{\shm}(C)} & {\ha{M}(C\cap\dist(U))=M_{C\cap\dist(U)}}
	\arrow["{i^*_C}", two heads, from=1-1, to=1-2]
	\arrow["{incl.}"', hook, from=2-1, to=1-1]
	\arrow["\exists", dashed, from=2-1, to=2-2]
	\arrow["{incl.}"', hook, from=2-2, to=1-2]
\end{tikzcd}\]

\begin{de}\label{msecdef}
    Let $X$ be a scheme locally of finite type, $C\subset\dist(X)$ be a subcoalgebra. We define We define $\ha{\shm}(C)\subset (C\sq\dist(\shm))^*$ to be the intersection
    $$\ha{\shm}(C)=\bigcap_{U}i^{*-1}_C(M_{C\cap\dist(U)})\subset(C\sq\dist(\shm))^*$$
    where $U$ runs over all the affine open subsets of $X$ and $M$ is an finitely generated $\O_X(U)$-module such that $\shm|_U\simeq\til{M}$. A $C^*$-module homomorphism $i_C$ is the natural map 
$$(C\cap\dist(U))\sq\dist(\til{M})\subset C\sq\dist(\shm).$$
    
    Note that $\ha{\shm}(C)$ is a module over $C^*$ and in particular a module over $\Sec_X(C)$.
\end{de}

In other words, the subspace $\ha{\shm}(C)\subset (C\sq\dist(\shm))^*$ is the   $\Sec_X(C)$-submodule formed by elements $\phi\in (C\sq \dist(\shm))^*$ satisfying the following condition:
\begin{itemize}
\item[(M)] For every affine open subscheme $U\subset X$, there exists $I\subset A$ such that $Z^\circ(I)\cap C\cap \dist(U)=0$ and $a_{C\cap\dist(U)}(I)\phi\subset m_{C\cap\dist(U)}(M)$ where $A=\O_X(U)$ and $\shm|_U=\til{M}$ for some finitely generated $\O_X(U)$-module $M$.
\end{itemize}

If $\dist(X)\supset C_1\supset C_2$ are subcoalgebras, then the inclusion $i:C_2\hookrightarrow C_1$ induces a homomorphism 
$$i\sq \id_{\shm}
:C_2\sq\dist(\shm)\rightarrow C_1\sq\dist(\shm)$$ 
of $C_1$-comodules. The dual map 
$$(i\sq \id_{\shm})^*:(C_1\sq\dist(\shm))^*\rightarrow (C_2\sq\dist(\shm))^*$$ 
restricts to a module homomorphism $\ha{\shm}(C_1)\rightarrow\ha{\shm}(C_2)$. In fact, let $U\subset X$ be an affine open subset and let $\shm|_U\simeq\til{M}$ for some finitely generated $\O_X(U)$-module $M$. Then the following diagram commutes:

\[\begin{tikzcd}
	{(C_1\sq\dist(\shm))^*} & {((C_1\cap\dist(U))\sq\dist(\til{M}))^*} \\
	{i_{C_1}^{*-1}(M_{C_1\cap\dist(U)})} & {\ha{M}(C_1\cap\dist(U))} & {(C_2\sq\dist(\shm))^*} & {((C_2\cap\dist(U))\sq\dist(\til{M}))^*} \\
	&& {i_{C_2}^{*-1}(M_{C_2\cap\dist(U)})} & {\ha{M}(C_2\cap\dist(U))}
	\arrow["i^*_{C_1}", from=1-1, to=1-2]
	\arrow["(i\sq\id_\shm)^*"{pos=0.8}, two heads, from=1-1, to=2-3]
	\arrow["(i'\sq\id_\shm)^*"{pos=0.8}, two heads, from=1-2, to=2-4]
	\arrow["incl."{pos=0.8}, hook, from=2-1, to=1-1]
	\arrow[from=2-1, to=2-2]
	\arrow[dashed, from=2-1, to=3-3]
	\arrow["incl."'{pos=0.8}, hook, from=2-2, to=1-2]
	\arrow[from=2-2, to=3-4]
	\arrow["i^*_{C_2}"', from=2-3, to=2-4]
	\arrow["incl."'{pos=0.8}, hook, from=3-3, to=2-3]
	\arrow[from=3-3, to=3-4]
	\arrow["incl.", hook, from=3-4, to=2-4]
\end{tikzcd}\]
where $i:C_2\hookrightarrow C_1$ and $i':C_2\cap\dist(U)\hookrightarrow C_1\cap\dist(U)$ are inclusions. We have seen in the previous subsection that the linear map
$$(i'\sq\id_\shm)^*:((C_1\cap\dist(U))\sq\dist(\til{M}))^*\rightarrow ((C_2\cap\dist(U))\sq\dist(\til{M}))^*$$ 
for restricts to the linear map
$\ha{M}(C_1\cap\dist(U))\rightarrow\ha{M}(C_2\cap\dist(U))$.  Therefore $(C_1\sq\dist(U))^*\rightarrow (C_2\sq\dist(U))^*$ restricts to 
$i_C^{*-1}(M_{C_1\cap\dist(U)})\rightarrow i_C^{*-1}(M_{C_2\cap\dist(U)})$.
Since $f$ is arbitrary, $C_1^*\rightarrow C_2^*$ restricts to a linear map
$$\ha{\shm}(C_1)=\bigcap_{U}i^{*-1}_{C_1}(M_{C\cap\dist(U)})\rightarrow\bigcap_{U}i^{*-1}_{C_2}(M_{C_2\cap\dist(U)})=\ha{\shm}(C_2)$$
Thus the construction of $\ha{\shm}(C)$ for every subcoalgebra $C\subset\dist(X)$ defines a functor $\ha{\shm}:P^{op}_{\dist(X)}\rightarrow \vect$.

\begin{de}
    Let $X$ be a scheme locally of finite type and $\shm$ be a coherent sheaf over $X$.
    Then the functor
$\ha{\shm}:P^{op}_{\dist(X)}\rightarrow\vect$ is a module over the ringed coalgebra $\dist(X)$. It will be called the module over $\dist(X)$ associated to $\shm$.
\end{de}

\begin{cor}
    Let $X$ be a scheme locally of finite type. The assignment that sends a coherent sheaf $\shm$ over $X$ to the module $\ha{\shm}$ over the ringed coalgebra $\emb(X)$ defines a functor $\ha{(-)}:\coh_X\rightarrow{}_{\emb(X)}\modu$.
\end{cor}
\begin{proof}
    Let $\shm_1$ and $\shm_2$ be coherent sheaves on $X$ and $f:\shm_1\rightarrow\shm_2$ be a morphism of sheaves. It induces a comodule morphism $\dist(\shm_2)\rightarrow\dist(\shm_1)$. For a subcoalgebra $C\subset \dist(X)$, one can show by an argument similar to \ref{tfunc} that the natural homomorphism $(C\sq\dist(\shm_1))^*\rightarrow(C\sq\dist(\shm_2))^*$ 
    of $C^*$-modules restricts to a map $\ha{\shm_1}(C)\rightarrow\ha{\shm_2}(C)$.
\end{proof}

Our next goal is to show that $\ha{(-)}$ is faithful under the condition that $X$ is separated. We first prove the following:

\begin{prop}\label{modres}
Let $X$ be a separated scheme locally of finite type and $U\subset X$ be an affine open subset. Then the following diagram commutes up to natural isomorphism:
\[\begin{tikzcd}
	{\coh_X} & {{}_{\emb(X)}\modu} \\
	{\coh_U} & {{}_{\emb(U)}\modu}
	\arrow["{\ha{(-)}}", from=1-1, to=1-2]
	\arrow["{res.}"', from=1-1, to=2-1]
	\arrow["{res.}"',from=1-2, to=2-2]
	\arrow["{\ha{(-)}}"', from=2-1, to=2-2]
\end{tikzcd}\]
\end{prop}
\begin{proof}
Let $\shm$ be a coherent sheaf on $X$, $U\subset X$ be an open subset and $C\subset\dist(U)$ be a subcoalgebra. We first check that $\ha{\shm|_U}(C)=\ha{\shm}(C)$. Note that $\dist(\shm)=\dist(\shm|_U)\oplus D$ and $\dist(U)\sq D=0$ for subcomodules $$D=\bigoplus_{x\in\lv X\backslash U\rv}\shm_x^\circ\subset\dist(\shm).$$ 
Since the cotensor product preserves direct sums, we have $C\sq\dist(\shm|_U)=C\sq\dist(\shm)$ and hence
$$(C\sq\dist(\shm|_U))^*=(C\sq\dist(\shm))^*.$$
We show that $\ha{\shm|_U}(C)=\ha{\shm}(C)$ as subsets of those dual spaces. By definition $\ha{\shm}(C)$ contains $\ha{\shm|_U}(C)$. Let $\phi\in \ha{\shm|_U}(C)$ satisfying the condition (M) for any affine open subset of $U$ and let $V\subset X$ be an affine open subset. By the separation of $X$, the intersection $U\cap V\subset X$ is affine. Also, $C\subset \dist(U)$ implies that
$$(C\cap\dist(V))\sq\dist(\shm|_V)=(C\cap\dist(U\cap V))\sq\dist(\shm|_{U\cap V}).$$ 
Since $U\cap V\subset U$ is affine open subset and $\phi$ satisfies (M), it is an element of $\ha{\shm}(C)$.

Next, let $f:\shm_1\rightarrow\shm_2$ be a morphism of coherent sheaves on $X$ and $f|_U:\shm_1|_U\rightarrow\shm_2|_U$  the morphism between the restrictions of the sheaves to $U$ and $C\subset\dist(U)$ a subcoalgebra. We need to show that the $C$-component of $\ha{f}$ and that of $\ha{f|_U}$ agree. However, both of them are given as restrictions and corestrictions of the same morphisms $(C\sq\dist(\shm_1))^*\rightarrow(C\sq\dist(\shm_2))^*$ of $C^*$-modules, i.e., the following diagram commutes:
\[\begin{tikzcd}
	{(C\sq\dist(\shm_1))^*} & {(C\sq\dist(\shm_2))^*} \\
	{(C\sq\dist(\shm_1|_U))^*} & {(C\sq\dist(\shm_2|_U))^*}
	\arrow["{f'}", from=1-1, to=1-2]
	\arrow["id", from=1-1, to=2-1]
	\arrow["id", from=1-2, to=2-2]
	\arrow["{f|_U'}"', from=2-1, to=2-2]
\end{tikzcd}\]
Here
$$f'=(id_C\sq\dist(f))^*, f|_U'=(id_C\sq\dist(f|_U))^*.$$
The proof is complete.
\end{proof}

\begin{prop}\label{minj}
Let $X$ be a separated scheme locally of finite type. The functor $\ha{(-)}:\coh_X\rightarrow{}_{\dist(X)}\modu$ is faithful, i.e., the map
$$\ha{(-)}:\coh_X(\shm_1, \shm_2)\rightarrow{}_{\emb(X)}\modu(\ha{\shm_1}, \ha{\shm_2})$$
is injective. 
\end{prop}
\begin{proof}
Let $f, g:\shm_1\rightarrow \shm_2$ be morphisms of $\coh_X$ such that $\ha{f}=\ha{g}$. It suffices to show that for every $x\in X$, $f$ and $g$ agree on some open neighborhood of $x$. Let $x\in X$ and $U$ be an affine open neighborhood of $x$. Then $\shm_1|_U\simeq\til{M}_1$ and $\shm_2|_U\simeq\til{M}_2$ for some finitely generated modules $M_1$ and $M_2$ over $\O_{X}(U)$. Let $f|_U, g|_U:\til{M}_1\rightarrow\til{M}_2$ be the induced morphisms between the sheaves restricted to $U$. Then by the preceding lemma, we have
$$\ha{f|_U}=\ha{f}|_{\dist(U)}=\ha{g}|_{\dist(U)}=\ha{g|_U}.$$
Thus $f|_U=g|_U$. Since the affine open subset $U\subset X$ is arbitrary we obtain $f=g$.
\end{proof}

%The $U$-components $f_U, g_U:M_1\rightarrow M_2$ of $f$ and $g$ induce $\ha{f_U}, \ha{g_U}:\ha{M}_1\rightarrow\ha{M}_2$. However, we have where $\rho$ stands for a restriction. 

\begin{lem}\label{afcoh}
Let $A$ be a commutative finitely generated algebra and $X=\spc(A)$. For every finitely generated $A$-module $M$, 
$$\dist(\til{M})(C)=\ha{M}(C).$$ 
In particular, $\dist(\til{M})$ and $\ha{M}$ are isomorphic as modules over the ringed coalgebra $A^\circ$.
\end{lem}
\begin{proof}
Applying condition (M) to the identity morphism on $X=\spc(A)$, we have $\dist(\til{M})(C)\subset M_C$. On the other hand, every element in $M_C$ satisfies condition (M) so that $M_C\subset\dist(\til{M})(C)$.
\end{proof}

\subsection{Modules over $\dist{(X)}$ when $k$ is algebraically closed}\label{mod3}
In this subsection, \textbf{we assume $k$ to be algebraically closed}. Let $X$ be a scheme locally of finite type over $k$. We have associated to each coherent sheaf $\shm$ on $X$ a module $\ha{\shm}$ over the ringed coalgebra $\dist(X)$ and have seen that the functor defined by this assignment is faithful if $X$ is separated. We will show that the functor $\ha{(-)}$ is full if $X$ is a separated scheme locally of finite type. Again, since $k$ is algebraically closed, the subset $\lv X\rv\subset X$ is the set of closed points of $X$.

Recall that for a scheme $X$ locally of finite type and a coherent sheaf $\shm$ on $X$, we have
$$\dist(\shm)=\bigoplus_{x\in\lv X\rv}\shm^\circ_{X, x}$$
so that
$$\dist(\shm)^*=\prod_{x\in\lv X\rv}\shm_{X, x}^{\circ*}=\prod_{x\in\lv X\rv}\hat{\shm}_{X, x}$$
where $\hat{\shm}_{x}$ is the completion of the left $\O_{X, x}$-module $\shm_{x}$ by the unique maximal ideal. Note that 
$$\dist(U)\sq\dist(\shm)=\bigoplus_{x\in\lv U\rv}\shm_x^\circ=\dist(\shm|_U)$$
for open subscheme $U\subset X$. 

\begin{lem}\label{modinj}
Let $X$ be a scheme locally of finite type, $U\subset X$ be an open subset and $\shm_X$ be an $\O_X$-module. The composition of the natural homomorphisms
$$\shm_X(U)\hookrightarrow\prod_{x\in U}\shm_{X, x}\twoheadrightarrow \prod_{x\in \lv U\rv}\shm_{X, x}$$
is injective.
\end{lem}
\begin{proof}
Let $f\in \shm_X(U)$ be a global section in the kernel of this homomorphism and $U=\bigcup_i U_i$ be an affine open covering. Then $\lv U\rv=\bigcup_i \lv U_i\rv$. We have the following commutative diagram:
\[\begin{tikzcd}
	{\shm_X(U)} && {\prod_{x\in \lv U\rv}\shm_{X, x}} \\
	\\
	{\shm_X(U_i)} && {\prod_{x\in \lv U_i\rv}\shm_{X, x}}
	\arrow[from=3-1, to=3-3]
	\arrow["{rest.}", from=1-1, to=3-1]
	\arrow[two heads, from=1-3, to=3-3]
	\arrow[from=1-1, to=1-3]
\end{tikzcd}\]
Here the horizontal arrow on the bottom is injective. In fact, the arrow can be described as the canonical homomorphism
$$M\rightarrow\prod_{\m\in\spm(R)}M_\m$$
where $R=\O_X(U_i)$ and $M=\shm(U_i)$. To show that this is injective, let $x\in M$ be an element of the kernel of the homomorphism. Then the annihilator $\ann_R(x)$ cannot be contained in any maximal ideal of $R$. Thus $1\in ann_R(x)$ and hence $x=0$. 

By the injectivity of the horizontal arrow on the bottom of the diagram, $f|_{U_i}=0$ as an element of $\shm_X(U_i)$. By the gluing property, we obtain $f=0$ as an element of $\shm_X(X)$. 
\end{proof}

Therefore we have a composite of inclusions
$$\shm_X(U)\hookrightarrow\prod_{x\in \lv U\rv}\shm_{X, x}\hookrightarrow\prod_{x\in \lv U\rv}\hat{\shm}_{X, x}=\dist(\shm|_U)^{*}$$
given by sending $f\in\shm_X(U)$ to a collection $\{f_x\}_{x\in \lv U\rv}$.
In this way, the $\O_X(U)$-module $\shm_X(U)$ can be seen as a $\O_X(U)$-submodule of $\dist(\shm|_U)^{*}$.

The following is an analogue of \ref{dsec}.

\begin{lem}\label{tfsec}
Let $A$ be a finitely generated algebra and $X=\spc(A)$. Then for every finitely generated $A$-module $M$, we have
$$\ha{M}(\dist(U))\simeq\til{M}(U)$$
for every open subscheme $U\subset X$.
\end{lem}
\begin{proof}
Note that $(\dist(U)\sq M^\circ)^*=\dist(\til{M}|_U)^*$. By definition, we have
$$M_{\dist(U)}=\{x\in (\dist(\til{M}|_U)^*|\exists I\in\F_{\dist(U)}\ a_{\dist(U)}(I)x\subset m^M_{\dist(U)}(M)\}$$
where $a_{\dist(U)}$ is the composite of natural algebra homomorphisms:
$$A\hookrightarrow A^{\circ*}\twoheadrightarrow \dist(U)^*=\prod_{x\in\lv U\rv}\hat{A}_x.$$
and $m^M_{\dist(U)}$ is the composite of natural module homomorphisms:
$$M\hookrightarrow M^{\circ*}\twoheadrightarrow \dist(\til{M}|_U)^*=\prod_{x\in\lv U\rv}\hat{M}_x.$$
%$$(\dist(U)\sq M^\circ)^*=\prod_{x\in\lv U\rv}\hat{M}_x.$$
Then the statement follows from an argument similar to the proof of \ref{dsec}.
\end{proof}

\begin{cor}\label{msec}
Let $X$ be a scheme locally of finite type, $Y\subset X$ an open subscheme of $X$ and $\shm$ a coherent sheaf on $X$. Then
$$\ha{\shm}(\dist(Y))\simeq\shm(Y).$$
\end{cor}
\begin{proof}
The scheme $Y$ is locally of finite type so $\shm(Y)$ can be viewed as a sub $\Gamma(Y)$-module of $\dist(\shm|_Y)^*$ by \ref{modinj}. Let $U\subset X$ be an affine open subscheme and $\shm(U)\rightarrow\dist(\shm|_U)$ be a corresponding $\O_X(U)$-module homomorphism. Then $\shm|_U\simeq\til{M}$ for some finitely generated module $M$ over $\O_X(U)$. Note that 
$$(\dist(Y)\sq\dist(\shm))^*=\dist(\shm|_Y)^*$$
and
$$((\dist(Y)\cap\dist(U))\sq\dist(\shm))^*=\dist(\shm|_{Y\cap U})^*\simeq\dist(\til{M}|_{Y\cap U})^*.$$
The dual $i_{\dist(Y)}^*:\dist(\shm|_Y)^*\rightarrow\dist(\til{M}|_{Y\cap U})^*$ of $i_{\dist(Y)}$ introduced before \ref{msecdef} is compatible with the restriction map 
$$\rho^Y_{Y\cap U}:\Gamma(Y)\rightarrow\Gamma(Y\cap U)\simeq\til{M}(Y\cap U)$$ 
induced by the morphism $f$ of schemes, i.e., the following commutes:
\[\begin{tikzcd}
	{\dist(\shm|_Y)^*} & {\dist(\til{M}|_{Y\cap U})^*} \\
	{\shm(Y)} & {\til{M}(Y\cap U)}
	\arrow["{i^*_{\dist(Y)}}", from=1-1, to=1-2]
	\arrow[hook, from=2-1, to=1-1]
	\arrow["{\rho^Y_{Y\cap U}}", from=2-1, to=2-2]
	\arrow[hook, from=2-2, to=1-2]
\end{tikzcd}\]
Since $Y\subset X$ is open, then $\ha{M}(\dist(Y\cap U_i))=\til{M}(Y\cap U)=\shm(Y\cap U)$ by \ref{tfsec}. Thus $i^*_{\dist(Y)}(\shm(Y))\subset \shm(Y\cap U)$ and hence $\shm(Y)$ is contained in $\ha{\shm}(\dist(Y))$.

%By definition the dual map 
%$$\dist(\shm|_Y)^*=(\dist(Y)\sq\dist(\shm))^*\rightarrow ((\dist(Y)\cap\dist(U))\sq\dist(\shm))^*=\dist(\shm|_{Y\cap U})^*=\dist(\til{M}|_{Y\cap U})^*$$
%restricts to a module homomorphism from $\ha{\shm}(\dist(Y))$ to $\ha{M}(\dist(Y\cap U_i))$ and to a restriction $\shm(Y)\rightarrow\shm(Y\cap U)$ of sections on $Y$ to those on $Y\cap U$. 
%=\ha{M}(Y\cap U)^*

We show the opposite containment. Take an element 
$\phi\in\ha{\shm}(\dist(Y))$ and an affine open covering $X=\bigcup_i U_i$. Let $\shm|_{U_i}\simeq \til{M}_i$ for some finitely generated 
$\O_X(U_i)$ module $M_i$. The dual map $j_i^*:\dist(\shm|_Y)^*\rightarrow\dist(\shm|_{Y\cap U_i})^*$ of the inclusion $j_i:\dist(\shm|_{Y\cap U_i})\subset\dist(\shm|_Y)$ restricts to a module homomorphism from $\ha{\shm}(\dist(Y))$ to $\ha{M}_i(\dist(Y\cap U_i))$. Since $Y$ is open, then $\ha{M}_i(\dist(Y\cap U_i))=\til{M}_i(Y\cap U_i)$ by \ref{dsec}. Then $\{j^*_i(\phi)\in\til{M}_i(Y\cap U_i)\}$ can be glued to obtain a unique element $\phi'\in\shm(Y)$. Since $\phi'\in\ha{\shm}(\dist(Y))$, then $j^*_i(\phi')=j^*_i(\phi)$ for all $i$ by construction.
Since the module homomorphism
$$\dist(\shm|_Y)^*\rightarrow \prod_i\dist(\shm|_{Y\cap U_i})^*$$
defined by $\psi\mapsto (j^*_i(\psi))_i$ is injective, we have $\phi=\phi'\in\Gamma(Y)$.
\end{proof}

%Let $\tfc_X$ denote the full subcategory of $\coh_X$ whose objects are torsion-free. 

By \ref{minj} and \ref{msec} we obtain the following:

\begin{cor}
Let $X$ be a separated scheme locally of finite type over an algebraically closed field $k$. The functor $\ha{(-)}:\coh_X\rightarrow{}_{\emb(X)}\modu$ is fully-faithful, i.e., the map
$$\ha{(-)}:\coh_X(\shm_1, \shm_2)\rightarrow{}_{\emb(X)}\modu(\ha{\shm_1}, \ha{\shm_2})$$
is bijective for any coherent sheaves $\shm_1$ and $\shm_2$.
\end{cor}
\begin{proof}
By \ref{minj}, it suffices to show the surjectivity.
For an open subset $U\subset X$, every morphism $g:\ha{\shm_1}\rightarrow\ha{\shm_2}$ induces a natural homomorphism $f_U:\shm_1(U)\simeq\ha{\shm_1}(\dist(U))\rightarrow\ha{\shm_2}(\dist(U))\simeq\shm_2(U)$ of $\O_X(U)$-modules. Here the isomorphisms follow from \ref{msec}.  Thus we may define a map ${}_{\emb(X)}\modu(\ha{\shm_1}, \ha{\shm_2})\rightarrow \tfc_X(\shm_1, \shm_2)$ by sending 
$g:\ha{\shm_1}\rightarrow\ha{\shm_2}$ 
to a natural transformation $f:=\{f_U\}_{U\subset X}:\shm_1\rightarrow\shm_2$. We show that $\ha{f}=g$. Let $U\subset X$ be an affine open subset. Then $\shm_1|_U\simeq\til{M}_1$ and $\shm_2|_U\simeq\til{M}_2$ for some finitely generated $\O_X(U)$-modules $M_1$ and $M_2$. We first show that the restrictions
$$\ha{f}|_{\dist(U)}, g|_{\dist(U)}:\ha{\shm}_1|_{\dist(U)}\rightarrow\ha{\shm}_2|_{\dist(U)}$$
agree. By \ref{modres}, $\ha{\shm}_1|_{\dist(U)}\simeq\ha{M}_1$ and $\ha{\shm}_2|_{\dist(U)}\simeq\ha{M}_2$. Lemma \ref{embmod} shows that  
the homomorphism $g':=\Gamma(g|_{\dist(U)}):M_1\rightarrow M_2$ of $\O_X(U)$-modules between the global sections of $\ha{\shm}_1|_{\dist(U)}$ and $\ha{\shm}_1|_{\dist(U)}$ satisfies $\ha{g'}=g|_{\dist(U)}$. However, the homomorphism $g'$ is nothing but $f_U$ by definition so that 
$$\ha{f}|_{\dist(U)}=\ha{f_U}=\ha{g'}=g|_{\dist(U)}.$$
Next, we show that the $C$-component $\ha{f}_C, g_C:\ha{\shm}_1(C)\rightarrow\ha{\shm}_2(C)$ of $\ha{f}$ and $g$ coincide. For every affine open subset $U\subset X$, the following commutative diagram of restrictions is commutative for $i=1, 2$:
\[\begin{tikzcd}
	{(C\sq\dist(\shm)_i)^*} & {(C\cap\dist(U))\sq\dist(\shm)_i)^*} \\
	{\ha{\shm}_i(C)} & {\ha{\shm}_i(C\cap\dist(U))}
	\arrow["res.", from=1-1, to=1-2]
	\arrow[hook, from=2-1, to=1-1]
	\arrow["res.", from=2-1, to=2-2]
	\arrow[hook, from=2-2, to=1-2]
\end{tikzcd}\]
%$$\ha{\shm}_i(C)\rightarrow\ha{\shm}_i(C\cap\dist(U)), i=1, 2.$$
Taking the product with respect to all the affine open subsets, one obtains
the following commutative diagram for $i=1, 2$:
\[\begin{tikzcd}
	{(C\sq\dist(\shm)_i)^*} & {\prod_U(C\cap\dist(U))\sq\dist(\shm)_i)^*} \\
	{\ha{\shm}_i(C)} & {\prod_U\ha{\shm}_i(C\cap\dist(U))}
	\arrow[from=1-1, to=1-2]
	\arrow[hook, from=2-1, to=1-1]
	\arrow[from=2-1, to=2-2]
	\arrow[hook, from=2-2, to=1-2]
\end{tikzcd}\]
Here the horizontal linear map at the top is injective. It is enough to show that 
$$C\sq\dist(\shm)_i=\bigcup_{U: \text{affine open}} (C\cap\dist(U))\sq\dist(\shm)_i).$$
Indeed, let $x\in C\sq\dist(\shm)_i$. Then the coaction $\rho(x)$ is in $C'\otimes\dist(\shm)_i$ for some subcoalgebra $C'\subset C$ of finite dimension. Since $X$ is quasi-separated, there exists an open affine subset $U\subset X$ that contains the finite set $pts(C')$ (\cite[\href{https://stacks.math.columbia.edu/tag/01ZU}{Tag 01ZU}]{SAO}). 
Thus the element $x$ is contained in $(C\cap\dist(U))\sq\dist(\shm)_i$. 

Therefore 
$$\ha{\shm}_i(C)\rightarrow\prod_{U:\text{affine\ open}}\ha{\shm}_i(C\cap\dist(U)),\ i=1, 2.$$
is also injective. Since $\ha{f}$ and $g$ are morphisms from $\ha{\shm}_1$ to $\ha{\shm}_2$, the following diagram commutes:
\[\begin{tikzcd}
	{\ha{\shm}_1(C)} & {\prod_U\ha{\shm}_1(C\cap\dist(U))} \\
	{\ha{\shm}_2(C)} & {\prod_U\ha{\shm}_2(C\cap\dist(U))}
	\arrow[hook, from=1-1, to=1-2]
	\arrow["{\ha{f}_C}"', from=1-1, to=2-1]
	\arrow["{g_C}", shift left=3, from=1-1, to=2-1]
	\arrow["h", from=1-2, to=2-2]
	\arrow[hook, from=2-1, to=2-2]
\end{tikzcd}\]
Here the vertical arrow $h$ on the right is the product
$$h:=\prod_U\ha{f}_{C\cap\dist(U)}=\prod_Ug_{C\cap\dist(U)}$$
where $U$ runs over all the open affine subsets in $X$. By the injectivity of the horizontal arrow at the bottom, we obtain $\ha{f}_C=g_C$. This completes the proof.
\end{proof}

\bibliographystyle{plain}
\bibliography{rcov}
\end{document}